\documentclass[12pt]{article}
\usepackage{CJK,bbm,amsfonts,mathrsfs,amsthm,amssymb,
amsmath, indentfirst,color} 

\def\dy{\displaystyle}
\def\e{\epsilon}
\def\ef#1 {$(\ref{#1})$}
\def\dy{\displaystyle}
\def\e{\epsilon}
\def\v{\varepsilon}
\def\x{\xi}
\def\t{\theta}

\def\k{\kappa}

\def\a{\alpha}
\def\b{\beta}
\def\g{\gamma}
\def\d{\delta}
\def\l{\lambda}

\def\f{\frac}

\def\r{\rho}

\def\z{\zeta}

\def\di{\displaystyle}
\def\i{\infty}
       \newtheorem{lemma}{\bf Lemma}[section]
       \newtheorem{theorem}[lemma]{\bf Theorem}

       \newtheorem{remark}[lemma]{\bf Remark}

\newcommand{\bn}{\begin{eqnarray}}
\newcommand{\en}{\end{eqnarray}}

\begin{document}
\begin{CJK}{GB}{gbsn}
\CJKtilde
\title{\bf The limit to rarefaction wave with vacuum for
 1D compressible fluids with temperature-dependent viscosities}
\author
 {\bf  Mingjie Li\thanks{The research of  M. Li was supported by the NSFC
Grant No. 11201503. E-mail: lmjmath@gmail.com.},~~Teng Wang\thanks{E-mail:
tengwang@amss.ac.cn.},~~and ~  Yi Wang\thanks{The research of Y. Wang was supported by the NSFC
Grant No. 11171326 and the CAS Program for Cross $\&$ Cooperative Team of  the Science $\&$ Technology Innovation. E-mail:
wangyi@amss.ac.cn.}\\
{\small $^\ast$College of Science, Minzu University of China,}\\
{\small Beijing 100081,
P. R. China}\\
{\small $^{\dag}$ Institute of Applied Mathematics,
AMSS,}\\
{\small CAS, Beijing 100190, P. R. China}\\
{\small$^\ddag$ Institute of Applied Mathematics and NCMIS,}\\{\small AMSS, CAS, Beijing 100190, P. R. China}\\
}
\date{}

 \maketitle

\begin{abstract}\noindent
In this paper we study the zero dissipation limit of the one-dimensional
full compressible Navier-Stokes(CNS) equations with temperature-dependent viscosity
and heat-conduction coefficient. It is proved that given a rarefaction wave with one-side vacuum state to the full compressible Euler equations, we can construct a sequence of solutions to  the full CNS equations which converge to the above rarefaction wave with vacuum as the viscosity and the heat conduction coefficient tend to zero. Moreover, the uniform convergence rate is obtained. The main difficulty in our proof lies in the degeneracies of the density, the temperature and the temperature-dependent viscosities at the vacuum region in the zero dissipation limit.
\end{abstract}

\textbf{Keywords}:~~compressible Navier-Stokes equations,
 temperature-dependent viscosities, zero dissipation limit, rarefaction wave, vacuum.

\bigskip

\textbf{    }

 \section{Introduction and main result} \setcounter{equation}{0}
We consider the zero dissipation limit of the one-dimensional compressible Navier-Stokes equations with heat-conduction in Eulerian
coordinates which read
\begin{equation}
\label{ns}
\left\{\begin{array}{l} \rho_{t}+(\rho u)_x=0,\quad\quad\qquad\qquad\qquad\quad x\in \mathbf{R}=(-\infty,+\infty), t>0,\\[2mm]
(\rho u)_t +\big(\rho u^2+p\big)_x=
\big(\e\mu(\t) u_x\big)_x,\\[2mm]
\dy[\r(e+\frac{1}{2}u^2)]_t+[\r u(e+\frac{1}{2} u^2)+up]_x=
(\e\kappa(\t)\theta_x+\e\mu(\t) uu_x)_x,
\end{array}\right.
\end{equation}
where
$\rho(x,t)\geq0,~ u(x, t),~p(x, t),~ e(x, t)\geq 0 $ and $ \theta(x, t)\geq
0$ represent the mass density, velocity, pressure, internal energy  and absolute
temperature of the gas, respectively, and $\e\mu(\t)$ and $\e\k(\t)$ denote
the viscosity and heat-conduction coefficient, respectively, with $\e>0$ being positive constant and
\begin{align}\label{heat}
\mu(\t)=\mu_1\t^\a, \qquad{\rm and}\qquad\k(\t)=\k_1 \t^\a
\end{align}
for positive constants $\mu_1$, $\k_1$ and $\a> 0$. Without loss of generality, it is assumed that $\mu_1=\k_1=1.$
Here we consider the ideal polytropic gas, that is, the pressure $p$ and the internal energy $e$ are given respectively by
\begin{align}\label{sta}
p=R\rho\theta=A\rho^\gamma\exp(\frac{\gamma-1}{R}S),\quad
e=\frac{R \theta}{\gamma-1},
\end{align}
satisfying the second law of thermodynamics
  \bn\label{st}
  de=\theta dS +\frac{p}{\rho^2}d\rho.
  \en
In the state equations \eqref{sta} and \eqref{st},
 $S=S(x,t)$ denotes the entropy of the gas and $\gamma>1$ is the adiabatic exponent and $A, R$ are both positive constants. For simplicity, it is normalized that  $$A=R=\g-1.$$
The compressible Navier-Stokes equations \eqref{ns} with the temperature-dependent viscosities \eqref{heat} can be derived exactly from the Chapman-Enskog expansion for the Boltzmann equation with respect to the Knudsen number, one can refer to \cite{Ch-C}, \cite{KMN} for the details. Formally, as
$\e\rightarrow 0+$, the system \eqref{ns} tends to the corresponding
inviscid Euler equations
\begin{equation}\label{eul}
\left\{\begin{array}{l} \rho_{t}+(\rho u)_x=0,\\[2mm]
(\rho u)_t +\big(\rho u^2+p\big)_x=0,\\[2mm]
\dy[\r(e+\frac{1}{2} u^2)]_t+[\r u(e+\frac{1}{2} u^2)+up]_x=0.
\end{array}\right.
\end{equation}
The Euler system \eqref{eul} is a strictly hyperbolic one for
 $\r>0$ whose first and third characteristic fields are genuinely nonlinear and second characteristic field is linearly degenerate, that
 is, in the equivalent system
$$
\left(
\begin{array}{l}
\di \r\\
\di u \\
\di S
\end{array}
\right)_t + \left(
\begin{array}{ccc}
\di u&\quad \r&\quad 0\\
 \f{p_\r}{\r}&\quad u&\quad \f{p_S}{\r}\\
\di 0& \quad 0 &\quad u
\end{array}\right)\left(
\begin{array}{l}
\di \r\\
\di u\\
\di S
\end{array}\right)_x=0,
$$
the Jacobi matrix
$$
\left(
\begin{array}{ccc}
\di u&\quad \r&\quad 0\\
 \f{p_{\r}}{\r}&\quad u&\quad \f{p_S}{\r} \\
\di 0&\quad 0 &\quad u
\end{array}\right)
$$
has three distinct eigenvalues
$$
\l_1(\r,u,\t)=u-\sqrt{p_\r(\r,S)},~
\l_2(\r,u,\t)=u,~
\l_3(\r,u,\t)=u+\sqrt{p_\r(\r,S)}
$$
with corresponding right eigenvectors
$$
r_1(\r,u,S)=(-\r,\sqrt{p_\r(\r,S)},0)^t,~ r_2(\r,u,S)=(p_S,0,-p_\rho)^t,~
r_3(\r,u,S)=(\r,\sqrt{p_\r(\r,S)},0)^t,
$$
such that
$$
r_i(\r,u,S)\cdot \nabla_{(\r,u,S)}\l_i(\r,u,S)\neq 0,\quad i=1,3,
$$
and
$$
r_2(\r,u,S)\cdot \nabla_{(\r,u,S)}\l_2(\r,u,S)\equiv 0.
$$
Thus the  two $i$-Riemann invariants $(i=1,3)$ can be defined by (cf. \cite{smoller})
\begin{equation}\label{RI}
\Sigma_i^{(1)}=u+(-1)^{\f{i-1}{2}}\int^{\r}\f{\sqrt{p_z(z,S)}}{z}dz,\qquad
\Sigma_i^{(2)}=S,
\end{equation}
such that
$$
\nabla_{(\r,u,S)} \Sigma_i^{(j)}(\r,u,S)\cdot r_i(\r,u,S)\equiv0,\quad i=1,3,~ j=1,2.
$$

The study of the limiting process of viscous flows when the
viscosity tends to zero  is one of the important problems in the
theory of the compressible fluid.
 When the solution of the inviscid flow is smooth, the zero
dissipation limit  can be solved by classical scaling method.
However, the inviscid compressible flow contains singularities such
as shock, contact discontinuity and the vacuum in general. Therefore, determining how to justify the
zero dissipation limit to the Euler equations with basic wave
patterns and the vacuum states is a natural and difficult problem.

There have been many results on the zero dissipation limit of the
compressible fluid with basic wave patterns without vacuum. For the
system of the hyperbolic conservation laws with artificial viscosity
$$
u_t+f(u)_x=\v u_{xx},
$$
Goodman-Xin \cite{good xin} first verified the viscous limit for
piecewise smooth solutions separated by non-interacting shock waves
using a matched asymptotic expansion method. Later Yu \cite{Yu}
proved it for the corresponding hyperbolic conservation laws with
both shock and initial layers. In 2005, important progress made by
Bianchini-Bressan\cite{B-B} justifies  the vanishing viscosity limit
in BV space even though the problem is still unsolved for the
physical system such as the compressible Navier-Stokes equations.

For the isentropic compressible Navier-Stokes equations with constant viscosity where the conservation of energy in \eqref{ns} is neglected,
Hoff-Liu \cite{hoffliu} first proved the vanishing viscosity limit
for piecewise constant shock even with initial layer. Later Xin
\cite{xin93} obtained the zero dissipation limit for rarefaction
waves without vacuum for both rarefaction wave data and
well-prepared smooth data. Then Wang \cite{Wang-H} generalized the
result of Goodmann-Xin \cite{good xin} to the isentropic
Navier-Stokes equations.
Recently, Chen-Perepelitsa \cite{chen-P} proved the vanishing
viscosity to the compressible Euler equations for the isentropic compressible
Navier-Stokes equations \eqref{ns} with constant viscosity by compensated compactness method
for the case that the far field of the initial values of Euler
system  has no vacuums.

For the full Navier-Stokes
equations \eqref{ns} with constant viscosity,
there are also many results on the zero dissipation limit to the
corresponding full Euler system with basic wave patterns without
vacuum. We refer to Jiang-Ni-Sun \cite{jiang} and Xin-Zeng
\cite{Xin-Zeng} for the rarefaction wave, Wang \cite{Wang} for the
shock wave, Ma \cite{Ma} for the contact discontinuity, Huang-Wang-Yang \cite{Huang-Wang-Yang} and Huang-Jiang-Wang \cite{H-J-W} for the
superposition of two rarefaction waves and a contact discontinuity,  Huang-Wang-Yang \cite{Huang-Wang-Yang-1}
for the superposition of one shock and one rarefaction wave and  Zhang-Pan-Wang-Tan \cite{ZPWT} for the superposition of two shock waves with the initial layer. Recently, Huang-Wang-Wang-Yang \cite{HWWY} succeed in justifies the vanishing viscosity limit of compressible Navier-Stokes equations in the setting of Riemann solutions for the superposition of shock wave, rarefaction wave and contact discontinuity.

 It is well-known that the vacuum states are generic in inviscid compressible Euler equations \eqref{eul} since the vacuum states may occur in the Riemann solutions instantaneously as $t>0$ even if the initial Riemann data are non-vacuum states on both sides at $t=0$. Therefore, vacuum states are important physical states in gas dynamics and often yield degeneracies and certain singularities in the physical system, which cause some essential analytical difficulties. For example, the velocity can not even be defined in the vacuum region. In the setting of Riemann solutions, as pointed out by Liu-Smoller \cite{lius}, among the three elementary hyperbolic waves, i.e., shock and rarefaction waves and contact discontinuities, to the
one-dimensional isentropic compressible Euler equations \eqref{eul},
only the rarefaction wave can be connected to the vacuum states. There are some mathematical results on the time-asymptotic stability and the vanishing viscosity limit to the rarefaction wave with the vacuum. Perepelitsa \cite{P} consider the time-asymptotic
stability of solutions to 1-d isentropic compressible Navier-Stokes equations with fixed and constant viscosity
toward rarefaction waves connected to vacuum in
Lagrangian coordinate. Then Jiu-Wang-Xin \cite{JWX} study the large
time asymptotic behavior
 toward rarefaction wave with vacuum for solution to
  the one-dimensional  compressible
Navier-Stokes equations with density-dependent viscosity, which can be viewed as the compressible Navier-Stokes equations \eqref{ns} with fixed $\e=1$ in isentropic regime. Note that in isentropic regime, there is no energy equation $\eqref{ns}_3$ and the viscosity coefficient in the momentum equation $\eqref{ns}_2$ transfer to depend on the density.   More recently, Huang-Li-Wang \cite{F-L-W} justified the vanishing viscosity limit of one-dimensional isentropic compressible
Navier-Stokes equations with constant viscosity to the rarefaction
wave with one-side vacuum state to the corresponding compressible
Euler equations. However, for the full compressible Navier-Stokes equations \eqref{ns} with temperature-dependent viscosities, as far as we know, there is no any result on the zero dissipation limit to the rarefaction wave with the vacuum due to various difficulties mentioned below.

Now we give a description of the rarefaction wave connected to the
vacuum to the full compressible Euler equations \eqref{eul}; see also the
reference \cite{smoller}. For definiteness,
3-rarefaction wave will be considered. If we investigate the
compressible
 Euler system \eqref{eul} with the Riemann initial data
\begin{equation}\label{Riemann}
\left\{
 \begin{array}{rr} \r(0,x)=0,&x<0,\\
(\rho,u,\t)(0,x)=(\r_+,u_+,\t_+),&x>0,
\end{array}
\right.
\end{equation}
where the left side is the vacuum state and $\r_+>0, u_+,\t_+>0$ are
prescribed constants on the right state, then the Riemann problem
\eqref{eul}, \eqref{Riemann} admits a $3-$rarefaction wave connected
to the vacuum on the left hand side. By the fact that along the
$3-$rarefaction wave curve, $3-$Riemann invariant $\Sigma^{(i)}_3(\r,u,\t),~(i=1,2)$ defined in \eqref{RI}
keeps constant in $(x,t)$, one can get the velocity
$$
u_-=\Sigma_3^{(1)}(\r_+,u_+,\t_+)
$$
 being the speed of the gas coming into
the vacuum from the 3-rarefaction wave. On the other hand, in the vacuum region, the absolute temperature $\t$ also becomes zero due to the state equation \eqref{sta}, that is, $\t=\r^{\g-1}e^S$ and the fact that the entropy $S$ keeps constant along the 3-rarefaction wave. Correspondingly, the main difficulty of the present paper lies in how to deal with the degeneracies of the temperature-dependent viscosities in the vacuum region in the zero dissipation limit process.

As described above, the $3-$rarefaction wave
connecting the vacuum state $\r=0$ to $(\r_+,u_+,\t_+)$ is the self-similar
solution $(\r^{r_3},u^{r_3},\t^{r_3})(\x),~(\x=\f xt)$ of \eqref{eul} defined
by
\begin{equation}
\begin{array}{c}
\l_3(\r^{r_3}(\x),u^{r_3}(\x),\t^{r_3}(\x))=\left\{
\begin{array}{ll}
\di\r^{r_3}(\x)\equiv0, &\di{\rm if}~~\x<\l_3(0,u_-,0)=u_-,\\
\di \x, &\di {\rm if}~~u_-\leq\x\leq\l_3(\r_+,u_+,\t_+),\\
\di \l_3(\r_+,u_+,\t_+), &\di {\rm if}~~\x>\l_3(\r_+,u_+,\t_+),
\end{array} \right.
\end{array}\label{r2}
\end{equation}
\begin{equation}
\Sigma_3^{(1)}(\r^{r_3}(\x),u^{r_3}(\x),\t^{r_3}(\x))=\Sigma_3^{(1)}(0,u_-,0)=\Sigma_3^{(1)}(\r_+,u_+,\t_+),\label{r2+}
\end{equation}
and
\begin{equation}
\Sigma_3^{(2)}=S^{r_3}=S_+:=-(\g-1)\log\r_+ +\log\t_+.\label{S2+}
\end{equation}

Thus one can define the momentum $m^{r_3}$ and the total internal energy $n^{r_3}$ of 3-rarefaction wave by
\begin{equation}m^{r_3}(\x):=\left\{\begin{array}{ll}\r^{r_3}(\x) u^{r_3}(\x),~~~ &{\rm if}~~~\r^{r_3}>0,\\
0,~~~&{\rm if}~~~ \r^{r_3}=0,
\end{array}\right.\label{r2++}
\end{equation}
and
\begin{equation}n^{r_3}(\x):=\left\{\begin{array}{ll}\r^{r_3}(\x) \t^{r_3}(\x),~~~ &{\rm if}~~~\r^{r_3}>0,\\
0,~~~&{\rm if}~~~ \r^{r_3}=0.
\end{array}\right.\label{e3}
\end{equation}
In the present paper, we want to construct a sequence of global-in-time solutions
$(\rho^\e,m^\e:=\r^{\e} u^{\e},n^\e:=\rho^\e \t^\e)(x, t)$ to the compressible Navier-Stokes equations
\eqref{ns} with temperature-dependent viscosities, which converge to the 3-rarefaction wave
$(\rho^{r_3},m^{r_3},n^{r_3})(x/t)$ with vacuum defined above as $\e$ tends to zero. The
effects of initial layers will be ignored by choosing the
well-prepared initial data \eqref{inie} depending on the viscosity for the
Navier-Stokes equations.

 As mentioned before, the main novelty and difficulty here is determining how to control the degeneracies caused by the vacuum in the rarefaction wave. To overcome this difficulty, we first cut off the 3-rarefaction wave with vacuum along the rarefaction wave curve suitably ($\nu$ is the cut-off parameter and the details can be seen in Section 2), and use the fact that the viscosity $\e$ can control the degeneracies
caused by the vacuum in rarefaction waves by choosing suitably
$\nu=\nu(\e)$. In fact, we choose $\nu=\e^{a}|\ln\e|$ with $a$
defined in \eqref{a} in the present paper. The other observation
is that we can carry out the energy estimates under the a priori
assumptions \eqref{assump}-\eqref{assump2} such that the perturbation is suitably small in
$L^\infty(\mathbf{R})$ norm with some decay rate with respect to $\e$.

On the other hand, compared with the previous works \cite{F-L-W} for the isentropic Naiver-Stokes equations with the constant viscosity case, some new difficulties occur for the full Navier-Stokes equations \eqref{ns} with temperature-dependent viscosities considered in the present paper. Firstly, in order to overcome the difficulties caused by the non-isentropic regime, the relative entropy-entropy flux pair $(\eta,q)$ defined in \eqref{entropy} is used as in \cite{liuyyz}. Secondly, since in the vacuum region, the temperature becomes zero due to the
entropy keeps constant by the structure of rarefaction wave. Therefore, not only the density but also the temperature cause the degeneracies in the vacuum region and so the temperature-dependent viscosities do. Thus, to deal with the terms for the temperature-dependent viscosities, such as \eqref{linearmu}, becomes subtle. In fact, the derivative estimates of the perturbation of the density depends on the second order
derivative estimates of velocity with some degenerate coefficients (see \eqref{step22}), which is quite different from the constant viscosity case  in \cite{F-L-W} and \cite{lw}. By choosing the convergence
rate $a$ suitably as in \eqref{a} and then the parameters $\nu, \d$ as in \eqref{mu} closes the a priori estimates and yields the desired result.

Now we state our main result as follows.

\begin{theorem}\label{thm1}
Let $(\rho^{r_3},m^{r_3},n^{r_3})(x/t)$ be the 3-rarefaction wave with one-side vacuum state defined by
\eqref{r2}-\eqref{e3}. Then there exists a small positive constant
$\epsilon_0$ such that for any $\epsilon\in(0,\epsilon_0)$, we can
construct a family of global smooth solutions $(\rho^\e,m^\e=\rho^\e
u^\e, n^\e=\rho^\e
\t^\e)(x,t)$ to the compressible
Navier-Stokes equation (\ref{ns}) satisfying the following properties. \\
(1)\begin{equation*}\begin{array}{rl}
(\rho^\e-\rho^{r_3}, m^\e-m^{r_3}, n^\e-n^{r_3}
),(\rho^\e,m^\e,n^\e)_x &\dy \in C(0,+\infty;L^2(\mathbf{R})),\\
u^\e_{xx},\t^\e_{xx} &\dy \in L^2(0,+\infty;  L^2(\mathbf{R})).
\end{array}\end{equation*}
(2) As the viscosity $\epsilon\rightarrow 0$, $(\rho^\e, m^\e,n^\e)(x,t)$
converges to $(\rho^{r_3}, m^{r_3},n^{r_3})(x/t)$ pointwisely except the
original point $(0,0)$. Furthermore, for any given positive constant
$l$, there exists a constant $C_l>0$, independent of $\epsilon$,
such that
\begin{equation}\label{cr}
\begin{array}{ll}
\dy\sup_{t\geq l}\|
\rho^\e(\cdot,t)-\rho^{r_3}(\frac{\cdot}{t})\|_{L^\infty}\leq
C_l~\epsilon^{a} |\ln\epsilon |,\\
\di \sup_{t\geq l}\|
m^\e(\cdot,t)-m^{r_3}(\frac{\cdot}{t})\|_{L^\infty}\leq
C_l~\epsilon^{a} |\ln\epsilon |,\\
\di \sup_{t\geq l}\|
n^\e(\cdot,t)-n^{r_3}(\frac{\cdot}{t})\|_{L^\infty}\leq
C_l~\epsilon^{a} |\ln\epsilon |,
\end{array}
\end{equation}
with the positive constant $a$ given by
\begin{equation}\label{a}
\dy a=\f{1}{18\g+12\a(\g-1)}.
\end{equation}

 \end{theorem}

\begin{remark}
 From \eqref{cr} and \eqref{a}, one can see that the decay rate $\epsilon^{a} |\ln\epsilon|$
in Theorem \ref{thm1}  decreases monotonically with respect to $\a$,
which is consistent to the observation that the viscosity becomes weaker
as $\a$ becomes larger due to the vacuum in the rarefaction wave, although the convergence rates in \eqref{cr} and \eqref{a} may not be optimal.
\end{remark}

\begin{remark}
By using some ideas in the present paper, Li-Wang \cite{lw} generalize Huang-Li-Wang \cite{F-L-W}'s result to the non-isentropic compressible Navier-Stokes equations \eqref{ns} with constant viscosity.
\end{remark}

The rest of the paper is organized as follows. In section 2, we
construct a smooth 3-rarefaction wave profile which approximates the cut-off
rarefaction wave for the Euler equations based on the inviscid Burgers equation. Then the global-in-time solution to CNS \eqref{ns} is obtained around the smooth 3-rarefaction wave profile and finally the
proof the Theorem \ref{thm1} is shown in Section 3.

Throughout this paper, $H^k(\mathbf{R}), k= 0,1,2,. . . $,  denotes
the $l$-th order Sobolev space with its norm
$$
\|f\|_k=(\sum^k_{j=0}\|\partial^j_yf\|^2)^\frac{1}{2}, \quad {\rm
and}~\|\cdot\|:=\|\cdot\|_{L^2(dz)},
$$
while $L^2(dz)$ means the $L^2$ integral over $\mathbf{R}$ with
respect to the Lebesgue measure $dz$, and $z=x$ or $y$. For
simplicity, we also write $C$ as generic positive constants which
are independent of time $t$ or $\tau$ and viscosity $\e$ unless otherwise
stated.

\section{Rarefaction waves} \setcounter{equation}{0}
Since there is no exact rarefaction wave profile for the
Navier-Stokes equations \eqref{ns}, the following approximate
rarefaction wave profile satisfying the Euler equations was
motivated by Matsumura-Nishihara \cite{mn86}  and Xin \cite{xin93}.

Consider the Riemann problem for the inviscid Burgers equation:
\begin{align}\label{bur}
\left\{\begin{array}{ll}
w_t+ww_x=0,\\
w(x,0)=\left\{\begin{array}{ll}
w_-,&x<0,\\
w_+,&x>0.
\end{array}
\right.
\end{array}
\right.
\end{align}
If $w_-<w_+$, then the Riemann problem $(\ref {bur})$ admits a
rarefaction wave solution $w^r(x, t) = w^r(\f xt)$ given by
\begin{align}\label{abur}
w^r(\f xt)=\left\{\begin{array}{lr}
w_-,&\f xt\leq w_-,\\
\f xt,&w_-\leq \f xt\leq w_+,\\
 w_+,&\f xt\geq w_+.
\end{array}
\right.
\end{align}
Motivated by \cite{mn86, xin93},  the approximate rarefaction wave to the
compressible Navier-Stokes equations \eqref{ns} can be constructed
by the Burgers equation
\begin{eqnarray}\label{dbur}
\left\{
\begin{array}{l}
\di w_{t}+ww_{x}=0,\\
\di w( 0,x
)=w_\d(x)=w(\f{x}{\d})=\f{w_++w_-}{2}+\f{w_+-w_-}{2}\tanh\f{x}{\d},
\end{array}
\right.\label{(2.11)}
\end{eqnarray}
where $\d>0$ is a small parameter
 to be determined and the hyperbolic tangent function $\tanh x=\f{e^x-e^{-x}}{e^x+e^{-x}}$. In fact, we choose $\d=\e^a$ in (\ref{mu}) with $a$ given by \eqref{a}. Note that the solution $w^r_\d(t,x)$ of the
problem (\ref{(2.11)}) can be given explicitly by
\begin{equation}\label{b-s}
w^r_\d(t,x)=w_\d(x_0(t,x)),\qquad x=x_0(t,x)+w_\d(x_0(t,x))t.
\end{equation}
And $w^r_\d(t,x)$ has the following properties:

\begin{lemma}(\cite{xin93, F-L-W})
\label{appr} The problem~$(\ref{dbur})$ has a unique smooth global
solution $w_\delta^r(x,t)$ for each $\delta>0$ such that
\begin{itemize}
\item[(1)] $w_-<w_\delta^r(x,t)<w_+, \ \partial_x w^r_\delta(x,t)>0,$
 \ for  $x\in\mathbf{R}, \ t\geq 0,\ \delta>0.$
\item[(2)] The following estimates hold for all $\ t> 0,\ \delta>0$ and
p$\in[1,\infty]$:
\begin{align}
\|\partial_x w^r_\delta(\cdot,t)\|_{L^p}\leq&\dy C
(w_+-w_-)^{1/p}(\delta+t)^{-1+1/p}, \label{w1}
\end{align}
\vspace{-9mm}
\begin{align}\label{w2}
  \|\partial^2_x w^r_\delta(\cdot,t)\|_{L^p}\leq&\dy
C(\delta+t)^{-1}\delta^{-1+1/p},
\end{align}
\vspace{-9mm}
\begin{align}\label{der22}
|\frac{\partial^2 w^r_\delta(x,t)}{\partial
x^2}|\leq\frac{4}{\delta}\frac{\partial w^r_\delta(x,t)}{\partial
x}.
\end{align}
\item[(3)] There exist a constant $\delta_0\in (0,1)$ such that for
$\delta\in(0,\delta_0], t>0$,
$$
\| w^r_\delta(\cdot,t)-w^r(\f\cdot t)\|_{L^\infty}\leq \dy C\delta
t^{-1}\big[\ln (1+t)+|\ln\delta |\big].
$$
\end{itemize}
\end{lemma}
Note that Lemma
\ref{appr} is a little different from the one in \cite{xin93} as stated in \cite{F-L-W}. For the detailed proof of Lemma \ref{appr}, one can refer to \cite{xin93} and \cite{F-L-W} and we omit it here for brevity.

As mentioned in the introduction, we will first cut off the 3-rarefaction
wave with vacuum along the wave curve in order to overcome the
degeneracies caused by the vacuum. More precisely, for any constant $\nu>0$ to
be determined, we can get a state $(\r,u,\t)=(\nu, u_\nu, e^{\bar S}\nu^{\g-1})$ belonging to
the 3-rarefaction wave curve, where $\bar S=S_+=-(\g-1)\log\r_+ +\log\t_+$. From the fact that 3-Riemann invariant
$\Sigma_3^{(i)}(\r,u,\t),~(i=1,2)$  is constant along the 3-rarefaction wave curve,
$u_\nu$ can be computed explicitly by
$$u_\nu=\Sigma_3^{(1)}(\r_+,u_+,\t_+)+2\sqrt{\frac{\g}{\g-1}\nu^{\g-1}e^{S_+}}.$$ Now we get
a new cut-off 3-rarefaction wave $(\rho_\nu^{r_3},u_\nu^{r_3},\t_\nu^{r_3}
)(\x),~(\x=x/t)$ connecting the state $(\nu,u_\nu,e^{\bar S}\nu^{\g-1})$ to the state
$(\r_+,u_+,\t_+)$ which can be expressed explicitly by
\begin{align}\label{rar1}
\lambda_3(\rho_\nu^{r_3},u_\nu^{r_3},\t_\nu^{r_3})(\xi)=\left\{\begin{array}{ll}
\lambda_3(\nu, u_\nu, e^{\bar S}\nu^{\g-1}),&\xi< \lambda_3(\nu, u_\nu, e^{\bar S}\nu^{\g-1}),\\
\xi,&\lambda_3(\nu, u_\nu, e^{\bar S}\nu^{\g-1})\leq \xi\leq \lambda_3(\rho_+,u_+,\t_+),\\
\lambda_3(\rho_+,u_+,\t_+),&\xi>\lambda_3(\rho_+,u_+,\t_+)
\end{array}
\right.
\end{align}
and
\begin{eqnarray}
\label{rar3}
\Sigma_3^{(1)}(\rho_\nu^{r_3},u_\nu^{r_3},\t_\nu^{r_3})=\Sigma_3^{(1)}(\nu,
u_\nu, e^{\bar S}\nu^{\g-1})=\Sigma_3^{(1)}(\r_+,u_+,\t_+).
\end{eqnarray}
Correspondingly, we can define the momentum function and the total internal energy
$m^{r_3}_\nu:=\r^{r_3}_\nu u^{r_3}_\nu$ and $n^{r_3}_\nu:=\r^{r_3}_\nu \t^{r_3}_\nu$
 respectively. It is easy to show that the
cut-off 3-rarefaction wave $(\rho_\nu^{r_3},m_\nu^{r_3},n_\nu^{r_3} )(x/t)$
converges to the original 3-rarefaction wave with vacuum
$(\rho^{r_3},m^{r_3},n^{r_3} )(x/t)$ in sup-norm with the convergence rate
$\nu$ as $\nu$ tends to zero. More precisely, it holds that

\begin{lemma}\label{cut-off}
There exist a constant $\nu_0\in (0,1)$ such that for
$\nu\in(0,\nu_0], t>0$,
$$
\| (\rho_\nu^{r_3},m_\nu^{r_3},n_\nu^{r_3} )(\cdot/t)-(\rho^{r_3},m^{r_3},n^{r_3}
)(\cdot/t)\|_{L^\infty}\leq \dy C\nu,
$$
where the positive constant $C$ is independent of $\nu$.
\end{lemma}

Now the approximate rarefaction wave
$(\bar{\r}_{\nu,\d},\bar{u}_{\nu,\d},\bar{\t}_{\nu,\d})(x,t)$ of the cut-off
3-rarefaction wave $(\rho_\nu^{r_3},u_\nu^{r_3},\t_\nu^{r_3})(\f xt)$ to
compressible Euler equations $(\ref {eul})$ can be defined by
\begin{eqnarray}
\left\{
\begin{array}{l}
\di w_+=\l_3(\r_+,u_+,\t_+),\quad w_-=\l_3(\nu,u_\nu,e^{\bar S}\nu^{\g-1}), \\
\di w_\d^r(x,t)= \l_3(\bar{\r}_{\nu,\d},\bar{u}_{\nu,\d},\bar{\t}_{\nu,\d})(x,t),\\
\di
\Sigma_3^{(1)}(\bar{\r}_{\nu,\d},\bar{u}_{\nu,\d},\bar{\t}_{\nu,\d})(x,t)=\Sigma_3^{(1)}(\rho_+,u_+,\t_+)=
\Sigma_3^{(1)}(\nu,u_\nu,e^{\bar S}\nu^{\g-1}),
\end{array} \right.\label{au}
\end{eqnarray}
where $w_\delta^r(x,t)$ is the solution of Burger's equation $(\ref
{dbur})$ defined in \eqref{b-s}. From now on, the subscription of
$(\bar\r_{\d,\nu}, \bar u_{\d,\nu},\bar\t_{\d,\nu})(x,t)$ will be  abbreviated as
$(\bar\r,\bar u,\bar\t)(x,t)$ for simplicity. Then the approximate 3-rarefaction wave $(\bar \r,\bar u,\bar \t)$ defined above satisfies
\begin{equation}\label{ar}
\left\{\begin{array}{l} \bar\rho_{t}+(\bar\rho \bar u)_x=0,\\[2mm]
(\bar\rho \bar u)_t +\big(\bar\rho \bar u^2+\bar p\big)_x=0,\\[2mm]
\dy\Big[\bar\r(\bar e+\frac{\bar u^2}{2})\Big]_t+\Big[\bar\r\bar u(\bar e+\frac{\bar u^2}{2})+\bar u \bar p\Big]_x=0,
\end{array}\right.
\end{equation}
where
$$
\bar p=R \bar\rho\bar\t=A\bar\rho^\g \exp{(\f{\g-1}{R}\bar S)},\quad {\rm and}~~ \bar e=\frac{R}{\g-1}\bar\t.
$$
The properties of the approximate rarefaction wave $(\bar\r,\bar
u,\bar\t)$ is listed without proof in the following Lemma.

\begin{lemma}
\label{appu} The approximate cut-off 3-rarefaction wave $(\bar
\r,\bar u, \bar\t)$ defined in \eqref{au} satisfies the following
properties:
\begin{itemize}
\item[(i)] $\bar u_x(x,t)=\f2{\g+1}(w_\d^r)_x>0,$
 \ for  $x\in\mathbf{R}, \ t\geq 0,$\\[3mm]
$\bar\r_x=\frac{1}{\sqrt{\g(\g-1)e^{S_+}}}\bar\r^{\frac{3-\g}{2}}\bar u_x$, and
$\bar\r_{xx}=\frac{1}{\sqrt{\g(\g-1)e^{S_+}}}\bar\r^{\f{3-\g}{2}}\bar
u_{xx}+\f{3-\g}{2\g(\g-1)e^{S_+}}\bar\r^{2-\g}(\bar u_x)^2$,\\[3mm]
$\bar\t_x=\sqrt{\f{\g-1}{\g}}\bar\t^{\frac{1}{2}}\bar u_x$, and
$\bar\t_{xx}=\sqrt{\f{\g-1}{\g}}\bar\t^{\f{1}{2}}\bar
u_{xx}+\frac{\g-1}{2\g}(\bar u_x)^2$.
\item[(ii)] The following estimates hold for all $\ t> 0,\ \delta>0$ and
p$\in[1,\infty]$:
\begin{equation*}\begin{array}{l}
\| \bar u_x(\cdot,t)\|_{L^p}\leq
C (w_+-w_-)^{1/p}(\delta+t)^{-1+1/p},\\[4mm]
  \|\bar u_{xx}(\cdot,t)\|_{L^p}\leq
C(\delta+t)^{-1}\delta^{-1+1/p}.
\end{array}\end{equation*}
\item[(iii)] There exist a constant $\delta_0\in (0,1)$ such that for
$\delta\in(0,\delta_0], t>0$,
\begin{equation*}\begin{array}{ll}
\| (\bar\r-\rho^{r_3}_\nu, \bar
u-u^{r_3}_\nu, \bar\t-\t^{r_3}_\nu)(\cdot,t)\|_{L^\infty}\leq \dy C\delta t^{-1}\big[\ln
(1+t)+|\ln\delta |\big].
\end{array}\end{equation*}
\end{itemize}
\end{lemma}

The proof of Lemma \ref{appu} can be got similarly as in \cite{F-L-W} and will be omitted for brevity.

\section{Proof of Theorem~\ref{thm1}}\setcounter{equation}{0}

In order to prove Theorem $\ref{thm1}$, we construct the global smooth solution
$(\rho^\e, u^\e,\t^\e)$ as the perturbation around the approximate
rarefaction wave $(\bar\rho, \bar u, \bar\t)$.  Consider the Cauchy problem
 (1.1) with the smooth initial data
\begin{equation}
\begin{aligned}\label{inie}
(\rho^\e, u^\e,\t^\e)(x,t=0)=(\bar\r, \bar u,\bar\t)(x,0).
\end{aligned}
\end{equation}
Then we introduce the perturbation
\begin{equation}\label{pert}
(\phi, \psi,\zeta)(y,\tau)=(\rho^\e, u^\e,\t^\e)(x,t)- (\bar\r, \bar u,\bar\t)(x,t),
\end{equation}
where $y,\tau$ are the scaled variables as
\begin{equation}\label{tao}
y=\frac{x}{\epsilon},\quad  \tau=\frac{t}{\epsilon},
\end{equation}
and $(\rho^\e, u^\e,\t^\e)$ is assumed to be the solution to the problem
 $(\ref {ns})$. For the simplicity of the notation, we will omit the superscription of $(\rho^\e, u^\e,\t^\e)$
  as $(\rho, u,\t)$ from now on if there is no confusion of the notation. Substituting $(\ref{pert})$ and $(\ref{tao})$
into the system $(\ref {ns})$ and using the equations for $(\bar\r, \bar u,\bar\t)$, one can obtain
\begin{equation}\label{mass}
\left\{
\begin{array}{ll}
\dy \phi_\tau+\r\psi_y+u\phi_y=-f,
\\
\dy \r\psi_\tau+\r u\psi_y+(\gamma-1)(\t\phi_y+\r\zeta_y)- \mu(\bar{\theta})\psi_{yy}=-g+\mu(\bar{\theta})_{y}u_y+\big((\mu(\theta)-\mu(\bar{\theta}))u_y\big)_y,
\\
\dy \r\zeta_\tau+ \r u\zeta_y+(\gamma-1)\r\theta\psi_y-\kappa(\bar{\theta})\zeta_{yy}=
-h+\kappa(\bar{\theta})_{y}u_y+\big((\kappa(\theta)-\kappa(\bar{\theta}))\t_y\big)_y+\mu(\theta)u_y^{2},
\end{array}
\right.
\end{equation}
with the initial data
\begin{equation}\label{init}
(\phi, \psi, \zeta)(y,0)=0,
\end{equation}
where
\begin{equation}\label{g}
\left\{
\begin{array}{l}
f={\bar u}_y\phi+{\bar\r}_y\psi,
 \\
\dy g=-\mu(\bar \theta) \bar u_{yy}+\r\psi \bar u_y+
(\gamma-1)(\bar \rho_y\zeta-\frac{\bar\r_y\bar\t\phi}{\bar\r}),
\\[1mm]
\dy h=\rho\psi\bar\t_y+(\gamma-1)\r\zeta \bar u_y-\kappa(\bar \theta)\bar \theta_{yy}.
\end{array}\right.
\end{equation}

We seek a global-in-time solution
$(\phi,\psi,\zeta)$ to the problem $(\ref{mass})-(\ref{g})$. To this
end,  the solution space for $(\ref{mass})-(\ref{g})$  is defined by
\begin{align}
\chi(0,\tau_1(\e))=\Big\{(\phi,\psi,\zeta)\Big|&(\phi,\psi,\zeta)\in
C([0,\tau_1(\e)];H^1(\mathbf{R})),\nonumber\quad
\phi_y\in L^{2}(0,\tau_1(\e);L^2(\mathbf{R})),\nonumber\\
&\psi_y,\zeta_y\in L^{2}(0,\tau_1(\e);H^1(\mathbf{R})) \Big\}\nonumber
  \end{align}
with $0<\tau_1(\e)\leq+\infty$.

\begin{theorem}\label{thm31}
There exist positive constants $\e_1$ and $C$ independent of $\e$,
such that if $0<\e\leq\e_1$,~ then the problem $(\ref{mass})-(\ref{init})$ admits a unique
global-in-time solution $(\phi,\psi,\zeta)\in \chi(0,+\i)$ satisfying
\begin{equation}\label{main1+}\begin{array}{ll}
&\dy\sup_{\tau\in[0,+\i)}\int_{\mathbf{R}}
  \Big(\bar \r^{\g-2}\phi^2+ \bar\r\psi^2+ \bar\r^{2-\g}\zeta^2\Big)(\tau,y)dy\\[3mm]
 &\dy +\int^{+\i}_{0}\int_{\mathbf{R}} \Big[\bar u_y\Big(\bar\r^{\g-2}\phi^2+ \bar\r\psi^2+ \bar\r^{2-\g}\zeta^2\Big) +\bar \theta^{\alpha}\psi_y^2+ \bar \theta^{\alpha-1}\zeta_y^2\Big]dyd\tau \leq  \e^{\f13},
\end{array}\end{equation}
\begin{equation}\label{main2+}\begin{array}{ll}
&\dy\sup_{\tau\in[0,+\i)}\int_{\mathbf{R}}
  \frac{\bar\t^{2\a}}{\bar\r^3}\phi_y^2dy+
  \int^{+\i}_{0}\int_{\mathbf{R}}\frac{\bar\theta^{\alpha+1}}{\bar\r^2}\phi_y^2dyd\tau
\leq  \e^{\f13-3a\g}|\ln\e|^{-3\g},
\end{array}\end{equation}
and
\begin{equation}\label{main3+}
\begin{array}{ll}
&\dy\sup_{\tau\in[0,+\i)}\int_{\mathbf{R}}(\psi_y^2+\z_y^2)dy+ \int^{+\i}_{0}\int_{\mathbf{R}}\Big[\bar u_y(\psi_y^2+\z_y^2)
+\frac{\bar\t^{\a}}{\bar\r}(\psi_{yy}^2+\z_{yy}^2)\Big]dyd\tau \leq\ \e^{\f19},
\end{array}\end{equation}
where $a$ is given by \eqref{a}.
\end{theorem}

In what follows, the analysis is always carried out under the a priori
assumptions
\begin{align}
\sup_{\tau\in[0,\tau_1(\e)]} \|\phi(\cdot,\tau)\|_{L^\infty}\leq \e^{a},\quad
\sup_{\tau\in[0,\tau_1(\e)]} \|\zeta(\cdot,\tau)\|_{L^\infty}\leq \e^{(\gamma-1)a},\label{assump}\\
\sup_{\tau\in[0,\tau_1(\e)]} \|\psi(\cdot,\tau)\|_{L^\infty}\leq \e^{a},\quad
\sup_{\tau\in[0,\tau_1(\e)]} \|(\psi_y,\z_y)(\cdot,\tau)\|\leq 1,\label{assump2}
\end{align}
where $a$ is given by \eqref{a}, $[0,\tau_1(\e)]$ is the time interval in which the solution exists and $\tau_1(\e)$ may depend on $\e$.

Take
\begin{equation}\label{mu}
\nu=\e^{a}|\ln\e|, \qquad\d=\e^a,
\end{equation}
in the sequel. Then it follows that $\nu\geq C \e^{a}$ with $C\geq\max\{2, \big(2 e^{-\bar S}\big)^{\f{1}{\g-1}}\}$ if $\e\ll1$. Under
the a priori assumption \eqref{assump}, one can get
\begin{align}\label{rup}
\frac{\bar\r}{2}\leq \r\leq\frac{3\bar\r}{2},\quad {\rm and}\quad \frac{\bar\t}{2}\leq \t\leq\frac{3\bar\t}{2}.
\end{align}
In fact, if $\e\ll1$, then one has
$$
\r=\bar\r+\phi\geq \bar\r-\|\phi\|_{L^\infty}\geq\bar\r-\e^{a}\geq
\bar\r-\frac{1}{2}\nu\geq\frac{\bar\r}{2},
$$
and
$$
\r=\bar\r+\phi\leq \bar\r+\|\phi\|_{L^\infty}\leq\bar\r+\e^{a}\leq
\bar\r+\frac{1}{2}\nu\leq\frac{3\bar\r}{2}.
$$
Similarly, note that $\bar\t=\bar\r^{\g-1} e^{\bar S}\geq\nu^{\g-1}e^{\bar S}$ by the definition of the rarefaction wave profile defined in \eqref{au}, it holds that
$$
\t=\bar\t+\z\geq \bar\t-\|\z\|_{L^\infty}\geq\bar\t-\e^{a(\g-1)}\geq
\bar\t-\f{e^{\bar S}}{2}\nu^{\g-1}\geq
\bar\t-\f{\bar\t}{2}=\f{\bar\t}{2} ,
$$
and
$$
\t=\bar\t+\z\leq \bar\t+\|\z\|_{L^\infty}\leq\bar\t+\e^{a(\g-1)}\leq
\bar\t+\f{e^{\bar S}}{2}\nu^{\g-1}\leq
\bar\t+\f{\bar\t}{2}=\f{3\bar\t}{2}.
$$

Since the proof for the local existence of the solution to
$(\ref{mass})-(\ref{g})$ is standard, we omit it for brevity. Note that in order to get the convergence rate of the local solution with respect to $\e$ as in \eqref{assump}, the local existence time interval, denoted by $[0,\tau_0]$ where $\tau_0$ may depend on $\e$, that is, $\tau_0=\tau_0(\e)$. The next step for the proof of Theorem \ref{thm31} is to extend the local solution to the global solution in $[0,\infty)$ for small but fixed viscosity coefficient and heat conduction coefficient $\e$. To do so, it is sufficient to show the following a priori estimates for fixed $\e$ with $0<\e\ll1$.

\begin{lemma}\label{len}
\textbf{ (A priori estimates)}\ \ Let $(\phi,\psi,\z)\in\chi(0,\tau_1(\e))$ be a solution to the problem
$(\ref{mass})-(\ref{g})$, where $\tau_1(\e)$ is the maximum existence time of the solution satisfying a priori assumptions \eqref{assump}-\eqref{assump2}. Then there exists a positive constant $\e_2$ such that if  $0<\e\leq\e_2$,~ then
\begin{equation}\label{main1}\begin{array}{ll}
&\dy\sup_{\tau\in[0,\tau_1(\e)]}\int_{\mathbf{R}}
  \Big(\bar \r^{\g-2}\phi^2+ \bar\r\psi^2+ \bar\r^{2-\g}\zeta^2\Big)(\tau,y)dy\\[3mm]
 &\dy +\int^{\tau_1(\e)}_{0}\int_{\mathbf{R}} \Big[\bar u_y\Big(\bar\r^{\g-2}\phi^2+ \bar\r\psi^2+ \bar\r^{2-\g}\zeta^2\Big) +\bar \theta^{\alpha}\psi_y^2+ \bar \theta^{\alpha-1}\zeta_y^2\Big]dyd\tau \leq  \e^{\f13},
\end{array}\end{equation}
\begin{equation}\label{main2}\begin{array}{ll}
&\dy\sup_{\tau\in[0,\tau_1(\e)]}\int_{\mathbf{R}}
  \frac{\bar\t^{2\a}}{\bar\r^3}\phi_y^2dy+
  \int^{\tau_1(\e)}_{0}\int_{\mathbf{R}}\frac{\bar\theta^{\alpha+1}}{\bar\r^2}\phi_y^2dyd\tau
\leq  \e^{\f13-3a\g}|\ln\e|^{-3\g},
\end{array}\end{equation}
and
\begin{equation}\label{main3}
\begin{array}{ll}
&\dy\sup_{\tau\in[0,\tau_1(\e)]}\int_{\mathbf{R}}(\psi_y^2+\z_y^2)dy+ \int^{\tau_1(\e)}_{0}\int_{\mathbf{R}}\Big[\bar u_y(\psi_y^2+\z_y^2)
+\frac{\bar\t^{\a}}{\bar\r}(\psi_{yy}^2+\z_{yy}^2)\Big]dyd\tau \leq\ \e^{\f19}.
\end{array}\end{equation}
\end{lemma}

\textbf{ Proof of Lemma \ref{len}:}\  \ The proof of Lemma \ref{len} consists of the following steps.\\
\underline{\it Step 1. }\quad First, as in \cite{liuyyz}, one can define the entropy-entropy flux pair $(\eta,q)$ as
\begin{equation}\label{entropy}
\left\{
\begin{array}{l}
\di\eta=-\bar \theta\{\r S-\bar\r\bar S-
\nabla_\textbf{X}(\r S)|_{\textbf{X}=\overline{\textbf{X}}}\cdot(\textbf{X}-\overline{\textbf{X}})\},\\[2mm]
\dy q=-\bar \theta\{\r uS-\bar\r \bar u\bar S- \nabla_\textbf{X}(\r
S)|_{\textbf{X}=\overline{\textbf{X}}}\cdot(\textbf{Y}-\overline{\textbf{Y}})\},
\end{array}\right.
\end{equation}
where
\begin{equation}
\begin{array}{ll}
\textbf{X}=\big(\r,\r u,\r(\frac{1}{2}|u|^2+\theta)\big)^t,\\
\textbf{Y}=\big(\r u,\r u^2+(\g-1)\r\theta,\r
u(\frac{1}{2}|u|^2+\g\theta)\big)^t.
\end{array}
\end{equation}
Since
\begin{equation}
\begin{array}{ll}
(\r S)_{\r}=S+\frac{|u|^2}{2\theta}-\g,\\
(\r S)_{m}=-\frac{u}{\theta},\\
(\r S)_{E}=\frac{1}{\theta},
\end{array}
\end{equation}
with $m=\r u$ and $E=\r(\t+\f{|u|^2}{2})$,
we get

\begin{equation}\label{g+}
\left\{
\begin{array}{l}
\di\eta=\r\theta- \bar\theta\r S+ \r[(\bar S-\g)\bar\theta+
\frac{1}{2}|u-\bar u|^2]
+(\g-1)\bar\r \bar \theta,\\[2mm]
~~~
=\dy(\gamma-1)\r\bar\theta\Phi(\frac{\bar\r}{\r})+\frac{\r\psi^2}{2}
+\r\bar\theta\Phi(\frac{\theta}{\bar\theta}),\\

\dy q=u\eta+(\g-1)(u-\bar u)(\r\theta-\bar\r\bar\theta),
\end{array}\right.
\end{equation}
where
\begin{equation}\label{qua}
\Phi(\eta)=\eta-\ln\eta-1.
\end{equation}
Direct computations yield
\begin{equation}\label{eq1}
 \eta_\tau+ q_y+ {\bar u}_yH
=\psi(\mu(\t)u_y)_y+\f{\z}{\t}\mu(\t)u_y^2+\f{\z}{\t}(\k(\t)\t_y)_y,
\end{equation}
where
\begin{equation}\label{ry1+}
\begin{array}{ll}
H=\dy\r(u-\bar u)^2+(\g-1)\r\bar\theta\Phi(\frac{\theta}{\bar\theta})+
(\g-1)^2\r\bar\theta\Phi(\frac{\bar\r}{\r})+
\sqrt{\frac{\g-1}{\g}}\bar\theta^{\frac{1}{2}}\r(u-\bar u)
\big((\g-1)\log\frac{\bar\r}{\r}+\log\frac{\theta}{\bar\theta}\big)\\[3mm]
\quad\di\geq (1-\b)\r(u-\bar
u)^2+(\g-1)\r\bar\t\big[\Phi(\f{\t}{\bar\t})+(\g-1)\Phi(\f{\bar\r}{\r})-\f{1}{4\b\g}
\big((\g-1)\log\f{\bar\r}{\r}+\log\f{\t}{\bar\t}\big)^2 \big]\\[3mm]
\end{array}\end{equation}
with $0<\b<1$ is the positive constant to be determined.
Let
$$
x_1=\f{\t}{\bar\t},\qquad x_2=\f{\bar\r}{\r},
$$
under the a priori assumptions \eqref{assump}, one has $x_1,x_2\sim 1$ as $\epsilon\rightarrow0.$
Consider the following function
$$
f^\b(x_1,x_2)=x_1-\log x_1-1+(\g-1)(x_2-\log x_2-1)-\f{1}{4\b\g}
\big((\g-1)\log x_2+\log x_1\big)^2.
$$
It is easy to check that
$$
f^\b(1,1)=f^\b_{x_1}(1,1)=f^\b_{x_2}(1,1)=0,
$$
and the Hessian matrix of $f^\b$ at point (1,1) is
$$
\nabla^2f^\b(1,1)=\left(
\begin{array}{cc}
\di 1-\f{1}{2\b\g} &\quad\di -\f{\g-1}{2\b\g}\\[3mm]
\di -\f{\g-1}{2\b\g} &\quad\di (\g-1)(1-\f{\g-1}{2\b\g})
\end{array}
\right).
$$
Thus the determinant of $\nabla^2f^\b(1,1)$ is
$$
\det\nabla^2f^\b(1,1)=(
\g-1)(1-\f1{2\b}).
$$
Take $\b=\f{3}{4}$, it is easy to see that $\nabla^2f^\b(1,1)$ is definitely positive near the point $(1,1)$.
So one has
$$
f^\b(x_1,x_2)\geq C_\b(x_1^2+x_2^2), \qquad {\rm as}~~x_1,x_2\sim1.
$$
Therefore, under the a priori assumptions \eqref{assump} and take $\b=\f34$ in \eqref{ry1+}, one can get
\begin{equation}
H\geq C~[\bar\r\psi^2+\frac{\bar\theta}{\bar\r}\phi^2+\frac{\bar\r}{\bar\theta}\zeta^2].
\label{H}
\end{equation}
Similarly, by the facts $\Phi(1)=\Phi^\prime(1)=0$, and $\Phi^{\prime\prime}(1)=1>0$, one has
$$
\eta\geq C~[\bar\r\psi^2+\frac{\bar\theta}{\bar\r}\phi^2+\frac{\bar\r}{\bar\theta}\zeta^2]
$$
under a a priori assumptions \eqref{assump}.
The right hand-side of \eqref{eq1} can be rewritten as
\begin{equation}
\begin{array}{ll}
\di\big[\mu(\t)\psi\psi_y+\f{\k(\t)}{\t}\di\z\z_y\big]_y-\mu(\t)\f{\bar\t}{\t}\di\psi_y^2-\k(\t)\f{\bar\t}{\t^2}\di\z_y^2
+\mu(\t)\psi\bar u_{yy}+\f{\k(\t)}{\t}\di\z\bar\t_{yy}\\[3mm]
\qquad +\di\f{\mu(\t)}{\t}\z(2\psi_y\bar u_y+\bar u_y^2)+\f{\k(\t)}{\t^2}\di\z\z_y\bar\t_y+\psi\mu(\t)_y\bar u_y
+\f{1}{\t}\di\k(\t)_y\z\bar\t_y.
\end{array}\label{3.31}
\end{equation}
Then integrating the equation \eqref{eq1} over ${\mathbf{R}}^1\times
[
0,\tau]$ and using \eqref{rup}, (\ref{eq1})-(\ref{3.31}) imply
\begin{equation}\begin{array}{ll}
&\dy\int_{\mathbf{R}}
  \Big(\bar\r^{\g-2}\phi^2+ \bar\r\psi^2+ \bar\r^{2-\g}\zeta^2\Big)(\tau,y)dy\\[3mm]
&\dy  +\int^\tau_{0}\int_{\mathbf{R}}\Big[\bar u_y\Big(\bar\r^{\g-2}\phi^2+ \bar\r\psi^2+ \bar\r^{2-\g}\zeta^2\Big) +\bar \theta^{\alpha}\psi_y^2+ \bar \theta^{\alpha-1}\zeta_y^2\Big]dyd\tau \\[4mm]
\leq & \dy C\int^\tau_{0}\int_{\mathbf{R}} |\bar \theta^{\alpha-1}\zeta||(\bar\theta_{yy},\bar u_y^2)|+
|\bar \theta^{\alpha}\psi\bar u_{yy}|+|\t^{\a-1/2}\psi\bar u_y^2|+|\bar u_y(\bar \theta^{\alpha-1}\zeta\psi_y,\bar\t^{\a-\f32}\zeta\zeta_y,\bar\t^{\a-1}\psi\zeta_y)| dyd\tau\\[4mm]
:= &\dy \sum_{i=1}^4I_i.
\end{array}\end{equation}

By Sobolev inequality and Lemma \ref{appu}, we can obtain
\begin{equation}\label{I1}\begin{array}{ll}
I_1&=\dy\int^\tau_{0}\int_{\mathbf{R}} |\bar \theta^{\alpha-1}\zeta|~|(\bar\theta_{yy},\bar u_y^2)|dyd\tau\dy\leq
C\nu^{1-\g}\int^\tau_{0}\|(\bar\theta_{yy},\bar u_y^2)\|_{L^1}\|\zeta\|^{1/2}\|\zeta_y\|^{1/2}d\tau\\[4mm]
& \dy \leq C\nu^{-\f{(4+\a)(\g-1)}{4}}\int^\tau_{0}\frac{1}{\tau+\delta/\epsilon}\|\zeta\|^{1/2}\|\sqrt{\bar\t^{\a-1}}\zeta_y\|^{1/2}
d\tau\\[4mm]
& \dy \leq \frac{1}{16}\int^\tau_{0} \|\sqrt{\bar\t^{\a-1}}\zeta_y\|^{2}d\tau+
C\nu^{-\f{(4+\a)(\g-1)}{3}}\int^\tau_{0}\|\zeta\|^{\frac{2}{3}}
\big(\frac{1}{\tau+\delta/\epsilon}\big)^{\frac{4}{3}}d\tau\\[4mm]
& \dy \leq \frac{1}{16}\int^\tau_{0} \|\sqrt{\bar\t^{\a-1}}\zeta_y\|^{2}d\tau+
C\nu^{-\f{(4+\a)(\g-1)+2}{3}} \sup_{0\leq\tau\leq\tau_1(\e)}\|\sqrt{\bar\r^{2-\g}}\z\|^{2/3}
\int^\tau_{0}\big(\frac{1}{\tau+\delta/\epsilon}\big)^{\frac{4}{3}}d\tau\\[4mm]
& \dy \leq \frac{1}{16}\int^\tau_{0} \|\sqrt{\bar\t^{\a-1}}\zeta_y\|^{2}d\tau+
\frac{1}{16}\sup_{0\leq\tau\leq\tau_1(\e)} \|\sqrt{\bar\r^{2-\g}}\zeta\|^2+
C\nu^{-\f{(4+\a)(\g-1)+2}{2}}(\frac{\e}{\delta})^{1/2}\\[4mm]
& \dy \leq \frac{1}{16}\int^\tau_{0} \|\sqrt{\bar\t^{\a-1}}\zeta_y\|^{2}d\tau+
\frac{1}{16}\sup_{0\leq\tau\leq\tau_1(\e)} \|\sqrt{\bar\r^{2-\g}}\zeta\|^2+\e^{\f13},
\end{array}\end{equation}
where we used the fact that
$$
\di C\nu^{-\f{(4+\a)(\g-1)+2}{2}}(\frac{\e}{\delta})^{1/2}
=C\e^{\f{1-a(4\g-1+\a(\g-1))}{2}}|\ln\e|^{-\f{(4+\a)(\g-1)+2}{2}}
\le \e^{\f13},\qquad {\rm
if}~~ \e\ll1.
$$
Similarly, it holds that
\begin{equation}\label{I2}\begin{array}{ll}
I_2&=\dy\int^\tau_{0}\int_{\mathbf{R}} |\bar \theta^{\alpha}\psi\bar u_{yy}|dyd\tau\dy\leq
C\int^\tau_{0}\|\bar u_{yy}\|_{L^1}\|\psi\|^{1/2}\|\psi_y\|^{1/2}d\tau\\[4mm]
& \dy \leq C\nu^{-\f{\a(\g-1)}{4}}\int^\tau_{0}\frac{1}{\tau+\delta/\epsilon}
\|\psi\|^{1/2}~\|\sqrt{\bar\t^{\a}}\psi_y\|^{1/2} d\tau\\[4mm]
& \dy \leq \frac{1}{16}\int^\tau_{0} \|\sqrt{\bar\t^{\a}}\psi_y\|^{2}d\tau
+C\nu^{-\f{\a(\g-1)}{3}}\int_0^\tau\big(\frac{1}{\tau+\delta/\epsilon}\big)^{\frac{4}{3}}
\|\psi\|^{2/3}d\tau\\[4mm]
& \dy \leq \frac{1}{16}\int^\tau_{0} \|\sqrt{\bar\t^{\a}}\psi_y\|^{2}d\tau
+C\nu^{-\f{\a(\g-1)+1}{3}}\sup_{0\leq\tau\leq\tau_1(\e)}\|\sqrt{\bar\r}\psi\|^{2/3}
\int_0^\tau\big(\frac{1}{\tau+\delta/\epsilon}\big)^{\frac{4}{3}}
d\tau\\[4mm]
& \dy \leq \frac{1}{16}\int^\tau_{0} \|\sqrt{\bar\t^{\a}}\psi_y\|^{2}d\tau+
\frac{1}{16}\sup_{0\leq\tau\leq\tau_1(\e)}\|\sqrt{\bar\r}\psi\|^2
+C\nu^{-\f{\a(\g-1)+1}{2}}(\frac{\e}{\delta})^{1/2}\\[4mm]
& \dy \leq \frac{1}{16}\int^\tau_{0} \|\sqrt{\bar\t^{\a}}\psi_y\|^{2}d\tau+
\frac{1}{16}\sup_{0\leq\tau\leq\tau_1(\e)}\|\sqrt{\bar\r}\psi\|^2+\e^{\f13},
\end{array}\end{equation}
where we used the fact that
$$
\di C\nu^{-\f{\a(\g-1)+1}{2}}(\frac{\e}{\delta})^{1/2}=C\e^{\f{1-a(\a(\g-1)+2)}{2}}|\ln\e|^{-\f{\a(\g-1)+1}{2}}\le \e^{\f13},\qquad {\rm
if}~~ \e\ll1.
$$
By using Lemma \ref{appu} yields
\begin{equation}\label{I3}\begin{array}{ll}
I_3&=\dy\int^\tau_{0}\int_{\mathbf{R}} |\bar \theta^{\alpha-1/2}\psi\bar u_{y}^2|dyd\tau\\[4mm]
& \dy \leq \frac{1}{16}\int^\tau_{0}\int_{\mathbf{R}}\bar u_y\bar\r\psi^2 dyd\tau
+ C\nu^{-\g}\int^\tau_{0}\|\bar u_y\|^3_{L^3} d\tau\\[4mm]
&\dy  \leq \frac{1}{16}\int^\tau_{0}\int_{\mathbf{R}}\bar u_y\bar\r\psi^2 dyd\tau
+ C\nu^{-\g}\int^\tau_{0}\big(\frac{1}{\tau+\delta/\epsilon}\big)^2 d\tau\\[4mm]
&\dy \leq \frac{1}{16}\int^\tau_{0}\int_{\mathbf{R}}\bar u_y\bar\r\psi^2 dyd\tau
+C\nu^{-\g}(\f{\v}{\d}).
\end{array}\end{equation}
And
\begin{equation}\label{I4}\begin{array}{ll}
I_{4}&\leq \dy C\int^\tau_{0}\int_{\mathbf{R}}\bar\r^{-\f12}\bar\t^{\f{\a-1}{2}}|\bar u_y| \big(|\sqrt{\bar\r^{2-\g}}\zeta||\bar\t^{\f\a2}\psi_y|,|\sqrt{\bar\r^{2-\g}}\zeta||\bar\t^{\f{\a-1}{2}}\z_y|,
|\sqrt{\bar\r}\psi||\bar\t^{\f{\a-1}{2}}\z_y|\big)dyd\tau \\[4mm]
& \dy \leq \frac{1}{16}\int^\tau_{0}\int_{\mathbf{R}} \big(\bar\theta^{\alpha}\psi_y^2 +\bar\theta^{\alpha-1}\z_y^2 \big)dyd\tau+
      C\frac{\e}{\nu^{\g}\delta}\int^\tau_{0}\int_{\mathbf{R}} \bar u_y\big(\bar\r^{2-\g}\zeta^2+\bar\r\psi^2\big) dyd\tau\\[4mm]
& \dy \leq \frac{1}{16}\int^\tau_{0}\int_{\mathbf{R}}
\big(\bar\theta^{\alpha}\psi_y^2 +\bar\theta^{\alpha-1}\z_y^2+\bar u_y\bar\r^{2-\g}\zeta^2+\bar u_y\bar\r\psi^2 \big)dyd\tau,     \end{array}\end{equation}
where we used the fact that
$$
\di C\f{\e}{\nu^{\g}\d}=C\e^{1-a-a\g}|\ln\e|^{-\g}\le C\e^{\f12}|\ln\e|^{-\g}\le\f{1}{16},\qquad {\rm
if}~~ \e\ll1.
$$
Combining (\ref{I1})-(\ref{I4}) and recalling \eqref{mu}
yield that
\begin{equation}\label{step1}\begin{array}{ll}
&\dy\sup_{\tau\in[0,\tau_1(\e)]}\int_{\mathbf{R}}
  \Big(\bar \r^{\g-2}\phi^2+ \bar\r\psi^2+ \bar\r^{2-\g}\zeta^2\Big)(\tau,y)dy\\[3mm]
 &\dy +\int^{\tau_1(\e)}_{0}\int_{\mathbf{R}} \Big[\bar u_y\Big(\bar\r^{\g-2}\phi^2+ \bar\r\psi^2+ \bar\r^{2-\g}\zeta^2\Big) +\bar \theta^{\alpha}\psi_y^2+ \bar \theta^{\alpha-1}\zeta_y^2\Big]dyd\tau \leq \e^{\f13}.
\end{array}\end{equation}

\underline{\it Step 2. }\quad Next we derive the estimation of $\phi_{y}$. Differentiating $(\ref {mass})_1$ with respect to $y$ and then multiplying the
resulted equation by $\mu^2(\bar\theta)\phi_y/\r^3$ imply that
\begin{equation}\label{m}\begin{array}{ll}
&\dy(\mu^2(\bar\theta)\frac{\phi_y^2}{2\r^3})_\tau +(\mu^2(\bar\theta)\frac{u\phi_y^2}{2\r^3})_y
+\mu^2(\bar\theta)\frac{\psi_{yy}\phi_y}{\r^2}\\[4mm]
=&\dy -\mu^2(\bar\theta)\frac{\phi_y}{\r^3}({\bar u}_{yy}\phi+{\bar\r}_{yy}\psi+2{\bar\r}_y\psi_y)+
\frac{\phi_y^2}{2\r^3}(\mu^2(\bar\theta)_\tau+\mu^2(\bar\theta)_y u).
\end{array}\end{equation}
Multiplying $(\ref{mass})_2$ by $\mu(\bar\theta)\phi_y/\r^2$ gives
\begin{equation}\label{v}\begin{array}{ll}
&\dy (\mu(\bar\theta)\frac{\psi\phi_y}{\r})_\tau-(\mu(\bar\theta)\frac{\psi\phi_\tau}{\r})_y-\mu(\bar\theta)\psi_y^2+
 (\gamma-1)\Big(\mu(\bar\theta)\frac{\theta\phi_y^2}{\r^2}+\mu(\bar\theta)\frac{\phi_y\zeta_y}{\r}\Big)\\[2mm]
&\dy -\mu^2(\bar\theta)\frac{\psi_{yy}\phi_y}{\r^2}+
 \mu(\bar\theta)\Big(-{\bar u}_y\frac{\psi_y\phi}{\r}
 +{\bar\r}_y{\bar u}_y\frac{\psi\phi}{\r^2}-\bar\r{\bar u}_{y}\frac{\psi\phi_y}{\r^2}\Big)+
 \mu(\bar\theta)g\frac{\phi_y}{\r^2}
=\mu(\bar\theta)\mu(\bar\theta)_y\frac{u_y\phi_y}{\r^2}\\[4mm]&\dy +
\mu(\bar\theta)_\tau\frac{\psi\phi_y}{\r}-\mu(\bar\theta)_y\frac{\psi\phi_{\tau}}{\r}-
\mu(\bar\theta)\frac{\psi^2}{\r^2}\bar\r_y(\phi_y+\bar\r_y)+((\mu(\theta)-\mu(\bar\theta))u_y)_y\mu(\bar\theta)\frac{\phi_y}{\r^2}.
 \end{array}\end{equation}
Combining $(\ref m)$ and $(\ref v)$ together,
then integrating the resulted equation over ${\mathbf{R}}^1\times
[0,\tau]$ imply
\begin{equation}
\begin{array}{ll}
&\dy\int_{\mathbf{R}}\Big(\mu^2(\bar\theta)\frac{\phi_y^2}{2\r^3}+\mu(\bar\theta)\frac{\psi\phi_y}{\r}\Big)dy +(\gamma-1)\int^\tau_{0}\int_{\mathbf{R}}
 \Big(\mu(\bar\theta)\frac{\theta\phi_y^2}{\r^2}+\mu(\bar\theta)\f{\phi_y\z_y}{\r}\Big)dyd\tau\\[5mm]
= &\dy\int^\tau_{0}\int_{\mathbf{R}}\Big\{ \mu(\bar\theta)\psi_y^2+
 \mu(\bar\theta)\Big({\bar u}_y\frac{\psi_y\phi}{\r}-{\bar\r}_y{\bar u}_y\frac{\psi\phi}{\r^2}+
 \bar\r{\bar u}_{y}\frac{\psi\phi_y}{\r^2}\Big) dyd\tau\\[5mm]
 &\dy
-\mu^2(\bar\theta)\frac{\phi_y}{\r^3}\Big({\bar
u}_{yy}\phi+{\bar\r}_{yy}\psi+2{\bar\r}_y\psi_y\Big)
-\mu(\bar\theta)g\frac{\phi_y}{\r^2}\Big\}dyd\tau\\[5mm]
&\dy +\int^\tau_{0}\int_{\mathbf{R}} \Big\{ \frac{\phi_y^2}{2\r^3}\Big(\mu^2(\bar\theta)_\tau+\mu^2(\bar\theta)_y u\Big)+
 \mu(\bar\theta)\mu(\bar\theta)_y\frac{u_y\phi_y}{\r^2}+
\mu(\bar\theta)_\tau\frac{\psi\phi_y}{\r}-\mu(\bar\theta)_y\frac{\psi\phi_{\tau}}{\r}\\[5mm]
&\dy -\mu(\bar\theta)\bar\r_y\frac{\psi^2}{\r^2}(\phi_y+\bar\r_y)+
 \big((\mu(\theta)-\mu(\bar\theta))u_y\big)_y\mu(\bar\theta)\frac{\phi_y}{\r^2}\Big\}dyd\tau.
\end{array}\label{3.48}
\end{equation}
 Combining (\ref{step1}) with \eqref{3.48}, and using the equation $\eqref{mass}_1$ and the fact that
\begin{equation}\label{linearmu}\begin{array}{ll}
&\dy\Big(\big(\mu(\theta)-\mu(\bar\theta)\big)u_y\Big)_y\\[4mm]
=&\dy\big((\mu(\theta)-\mu(\bar\theta)\big)(\psi_{yy}+\bar u_{yy})+
\alpha(\theta^{\alpha-1}\theta_y-\bar\theta^{\alpha-1}\bar\theta_y)(\psi_y+\bar u_y)\\[3mm]
=&\dy\big((\mu(\theta)-\mu(\bar\theta)\big)(\psi_{yy}+\bar u_{yy})+
\alpha\big[\theta^{\alpha-1}\zeta_y\psi_y+\theta^{\alpha-1}\zeta_y\bar u_y\\[3mm]
&\di\qquad\qquad\qquad\qquad\qquad\qquad +(\theta^{\alpha-1}-\bar\theta^{\alpha-1})\bar\theta_y\psi_y
+(\theta^{\alpha-1}-\bar\theta^{\alpha-1})\bar\theta_y\bar u_y\big],
\end{array}\end{equation}
it holds that
\begin{equation}\begin{array}{ll}
&\dy\int_{\mathbf{R}}
  \Big(\mu^2(\bar\theta)\frac{\phi_y^2}{\bar\r^3}+\bar\r^{\g-2}\phi^2+ \bar\r\psi^2+
  \bar\r^{2-\g}\zeta^2\Big)(\tau,y)dy\\[5mm]
+&\dy \int^\tau_{0}\int_{\mathbf{R}}\Big[\bar u_y\Big(\bar\r^{\g-2}\phi^2+ \bar\r\psi^2+
\bar\r^{2-\g}\zeta^2\Big) +\bar \theta^{\alpha}\psi_y^2+
\bar \theta^{\alpha-1}\zeta_y^2+\frac{\bar\theta^{\alpha+1}}{\bar\r^2}\phi_y^2\Big]dyd\tau\\[4mm]
\leq &\dy C |\int^\tau_{0}\int_{\mathbf{R}}\Big(\mu(\bar\t)\bar u_y\f{\psi_y\phi}{\r}+\mu(\bar\t)\bar u_y\bar\r\f{\psi\phi_y}{\r^2} -2\mu^2(\bar\t)\f{\phi_y}{\r^3}\bar\r_y\psi_y\Big)dyd\tau|\\[5mm]
&\dy +C|\int_0^\tau\int_{\mathbf{R}}
\mu^2(\bar\t)\f{\phi_y}{\r^3}\bar u_{yy}\phi+\mu^2(\bar\t)\f{\phi_y}{\r^3}\bar \r_{yy}\psi
+\mu(\t)\f{\phi_y}{\r^2}\big(\mu(\t)-\mu(\bar\t)\big)\bar u_{yy} dyd\tau|\\[5mm]
&\dy+C|\int_0^\tau\int_{\mathbf{R}}\mu(\bar\t)g\f{\phi_y}{\r^2} dyd\tau|
+C|\int_0^\tau\int_{\mathbf{R}}\f{\phi_y^2}{2\r^3}\big(\mu^2(\bar\t)_{\tau}+\mu^2(\bar\t)_y u\big) dyd\tau|\\[5mm]
&\dy +C|\int_0^\tau\int_{\mathbf{R}} \Big\{\mu(\bar\t)\bar\r_y\bar u_y\f{\psi\phi}{\r^2}+\mu(\bar\t)\mu(\bar\t)_y\f{\psi_y\phi_y}{\r^2}
-\mu(\bar\t)_y\f{\psi(\r\psi_y+u \phi_y+\bar u_y\phi+\bar\r_y\psi)}{\r}\\[5mm]
&\dy-\mu(\bar\t)\bar\r^2_y\f{\psi^2}{\r^2}
+\a\big(\t^{\a-1}\z_y\bar u_y+(\t^{\a-1}-\bar\t^{\a-1})\bar\t_y\psi_y \big)
\mu(\bar\t)\f{\phi_y}{\r^2}\Big\} ~dyd\tau|\\[5mm]
&\dy +C|\int_0^\tau\int_{\mathbf{R}}\Big(\mu(\bar\t)_{\tau}\f{\psi\phi_y}{\r}-\mu(\bar\t)\bar\r_y\f{\psi^2}{\r^2}\phi_y\Big) ~dyd\tau|
\\[5mm]
&\dy +C|\int_0^\tau\int_{\mathbf{R}}\Big(\mu(\bar\t)\mu(\bar\t)_y\f{\bar u_y\phi_y}{\r^2}
+\a(\t^{\a-1}-\bar\t^{\a-1})\bar\t_y\bar u_y\mu(\bar\t)\f{\phi_y}{\r^2}\Big) ~dyd\tau|
\\[5mm]
&\dy+C|\int_0^\tau\int_{\mathbf{R}}\mu(\bar\t)\f{\phi_y}{\r^2}\big(\mu(\t)-\mu(\bar\t) \big)\psi_{yy} dyd\tau|
+C|\int_0^\tau\int_{\mathbf{R}}\t^{\a-1}\z_y\psi_y\mu(\bar\t)\f{\phi_y}{\r^2} ~dyd\tau|
\\[5mm]
&\dy +C\e^{\f13}:=\sum_{i=1}^9J_i+C\e^{\f13}.
\end{array}\label{long}
\end{equation}
The terms on the right-hand side of \eqref{long} will be estimated one by one as follows. By Lemma 2.3,
\eqref{rup}, and Cauchy inequality, it holds that
\begin{equation}\label{J1}\begin{array}{ll}
J_1&\le\dy C|\int^\tau_{0}\int_{\mathbf{R}}\Big\{\bar\r^{\f{\a(\g-1)-\g}{2}}\bar u_y
\big(|\sqrt{\bar\t^{\a}}\psi_y||\sqrt{\bar\r^{\g-2}}\phi|,|\sqrt{\f{\bar\t^{\a+1}}{\bar\r^2}}\phi_y||\sqrt{\bar\r}\psi| \big) \\[3mm]
&\dy \quad+\bar\r^{\a(\g-1)-\g}\bar u_y |\sqrt{\f{\bar\t^{\a+1}}{\bar\r^2}}\phi_y| |\sqrt{\bar\t^{\a}}\psi_y|\Big\} dyd\tau|
\\[3mm]
&\dy \le \f{1}{16}\int_0^\tau\int_{\mathbf{R}}(\f{\bar\t^{\a+1}}{\bar\r^2}\phi_y^2+\bar\t^{\a}\psi^2_y)~dyd\tau
+C\f{\e}{\nu^{\g}\d}\int_0^\tau\int_{\mathbf{R}}\bar u_y(\bar\r^{\g-2}\phi^2+\bar\r\psi^2)~dyd\tau,
\end{array}\end{equation}
Recalling \eqref{der22} from Lemma 2.1 and the fact $(i)$ in Lemma 2.3, one can arrive at
\begin{equation}\begin{array}{ll}
|\bar\r_{xx}|\le C(\bar\r^{\f{3-\g}{2}}\di\f{\bar u_x}{\d}+\bar\r^{2-\g}\bar u^2_x).
\end{array}\end{equation}
Thus one has
\begin{equation}\label{J2}\begin{array}{ll}
J_2&\le\dy C|\int^\tau_{0}\int_{\mathbf{R}}\bar\r^{\f{3\a(\g-1)-3\g}{2}}\f{\e}{\d}\bar u_y
|\sqrt{\f{\bar\t^{\a+1}}{\bar\r^2}}\phi_y|(|\sqrt{\bar\r^{\g-2}}\phi|+|\sqrt{\bar\r}\psi|+|\sqrt{\bar\r^{2-\g}}\z|)dyd\tau | \\[3mm]
&\dy \le \f{1}{16}\int_0^\tau\int_{\mathbf{R}}\f{\bar\t^{\a+1}}{\bar\r^2}\phi_y^2~dyd\tau
+C\big(\f{\e}{\nu^{\g}\d}\big)^3\int_0^\tau\int_{\mathbf{R}}\bar u_y(\bar\r^{\g-2}\phi^2+\bar\r\psi^2+\bar\r^{2-\g}\z^2)~dyd\tau.
\end{array}\end{equation}
And
\begin{equation}\label{J3}\begin{array}{ll}
J_3&=\dy C|\int^\tau_{0}\int_{\mathbf{R}}\mu(\bar\t)g\f{\phi_y}{\r^2} dyd\tau|
\le  \f{1}{16}\int_0^\tau\int_{\mathbf{R}}\f{\bar\t^{\a+1}}{\bar\r^2}\phi_y^2~dyd\tau
+C\int_0^\tau\int_{\mathbf{R}}\bar\t^{\a-1}\bar\r^{-2}g^2 dyd\tau.
\end{array}\end{equation}
Recalling \eqref{g}, \eqref{rup} and $(i)$ in Lemma 2.3, one can get
\begin{equation}\begin{array}{ll}
|g|&\le\dy  C\big(\bar\t^{\a}|\bar u_{yy}|+|\bar\r\bar u_y\psi|+|\bar\r_y\z|+|\bar\r_y\bar\r^{\g-2}\phi| \big)\\[2mm]
&\dy\le C \big(\bar\t^{\a}|\bar u_{yy}|+\bar u_y(|\bar\r\psi|+|\bar\r^{\f{3-\g}{2}}\z|+|\bar\r^{\f{\g-1}{2}}\phi|) \big).
\end{array}\end{equation}
Thus the last term in \eqref{J3} can be estimated by
\begin{equation}\begin{array}{ll}
&\dy|\int_0^\tau\int_{\mathbf{R}}\bar\t^{\a-1}\bar\r^{-2}g^2 dyd\tau|\\[2mm]
&\dy\le C\nu^{-1-\g}\int_0^\tau\|\bar u_{yy}\|^2 d\tau+C\int_0^\tau\int_{\mathbf{R}}\bar\t^{\a-1}\bar\r^{-1}
\bar u^2_y(\bar\r^{\g-2}\phi^2+\bar\r\psi^2+\bar\r^{2-\g}\z^2) dyd\tau
\\[3mm]
&\dy\le C\nu^{-1-\g}\big(\f{\e}{\d}\big)^2+C\f{\e}{\nu^{\g}\d}\int_0^\tau\int_{\mathbf{R}}
\bar u_y(\bar\r^{\g-2}\phi^2+\bar\r\psi^2+\bar\r^{2-\g}\z^2) dyd\tau.
\end{array}\label{g2}\end{equation}
By Cauchy inequality, it holds that
\begin{equation}\label{J4}\begin{array}{ll}
J_4&\dy\le C|\int_0^\tau\int_{\mathbf{R}}\f{\bar\t^{\a+1}}{\bar\r^2}\phi_y^2~\bar\t^{\a-\f32}\bar\r^{-1}\bar u_y dyd\tau|
\le C\nu^{\f{1-3\g}{2}}\f{\e}{\d}\int_0^\tau\int_{\mathbf{R}}\f{\bar\t^{\a+1}}{\bar\r^2}\phi_y^2dyd\tau\\[3mm]
&\dy\le \f{1}{16} \int_0^\tau\int_{\mathbf{R}}\f{\bar\t^{\a+1}}{\bar\r^2}\phi_y^2dyd\tau
\qquad {\rm if}~~ \e\ll1.
\end{array}\end{equation}
Recalling  \eqref{rup} and $(i)$ in Lemma 2.3, one can get
\begin{equation}\label{J5}\begin{array}{ll}
J_5&\dy\le C|\int_0^\tau\int_{\mathbf{R}}\Big\{\bar\r^{\a(\g-1)-\g}\bar u_y|\sqrt{\f{\bar\t^{\a+1}}{\bar\r^2}}\phi_y|~|\sqrt{\bar\t^{\a}}\psi_y|+\bar\r^{\a(\g-1)-\g}\bar u_y|\sqrt{\f{\bar\t^{\a+1}}{\bar\r^2}}\phi_y|~|\sqrt{\bar\t^{\a-1}}\z_y|\\[3mm]
&\dy\quad+\bar\r^{\f{\a(\g-1)-\g}{2}}\bar u_y|\sqrt{\bar\r}\psi|~|\sqrt{\bar\t^{\a}}\psi_y|
+\bar\r^{\f{\a(\g-1)}{2}-\g}\bar u_y|\sqrt{\f{\bar\t^{\a+1}}{\bar\r^2}}\phi_y|~|\sqrt{\bar\r}\psi|\\[3mm]
&\dy \quad+\bar\r^{\a(\g-1)-\g}\bar u_y(|\sqrt{\bar u_y\bar\r}\psi|~|\sqrt{\bar u_y\bar\r^{\g-2}}\phi|+\bar u_y\bar\r\psi^2)\Big\} dyd\tau|\\[3mm]
&\dy\le (\f{1}{32}+C\f{\e}{\nu^{\g}\d})\int_0^\tau\int_{\mathbf{R}}\f{\bar\t^{\a+1}}{\bar\r^2}\phi^2_y dyd\tau
+C\f{\e}{\nu^{2\g}\d}\int_0^\tau\int_{\mathbf{R}}(
\bar\t^{\a}\psi^2_y+\bar\t^{\a-1}\z^2_y+\bar u_y\bar\r\psi^2+\bar u_y\bar\r^{\g-2}\phi^2) dyd\tau.
\end{array}\end{equation}
Similarly, $J_6$ can be estimated as
\begin{equation}\label{J6}\begin{array}{ll}
J_6&\dy\le C |\int_0^\tau\int_{\mathbf{R}}\bar\t^{\a/2-1}\bar\r^{-1/2}\bar u_y
|\sqrt{\f{\bar\t^{\a+1}}{\bar\r^2}}\phi_y|~|\sqrt{\bar\r}\psi| dyd\tau|\\[3mm]
&\dy\le \f{1}{16} \int_0^\tau\int_{\mathbf{R}}\f{\bar\t^{\a+1}}{\bar\r^2}\phi_y^2dyd\tau
+C\nu^{1-2\g}\f{\e}{\d} \int_0^\tau\int_{\mathbf{R}}\bar u_y\bar\r\psi^2 dyd\tau.
\end{array}\end{equation}
Recalling  \eqref{rup} and Lemma 2.3, one can get
\begin{equation}\label{J7}\begin{array}{ll}
J_7&\dy\le \f{1}{16} \int_0^\tau\int_{\mathbf{R}}\f{\bar\t^{\a+1}}{\bar\r^2}\phi_y^2dyd\tau
+C \int_0^\tau\int_{\mathbf{R}}\bar\t^{3\a-2}\bar\r^{-2}\bar u_y^4 dyd\tau\\[3mm]
&\dy\le \f{1}{16} \int_0^\tau\int_{\mathbf{R}}\f{\bar\t^{\a+1}}{\bar\r^2}\phi_y^2dyd\tau
+C \nu^{-2\g}\int_0^\tau\|\bar u_y\|^4_{L^4} d\tau\\[3mm]
&\dy\le \f{1}{16} \int_0^\tau\int_{\mathbf{R}}\f{\bar\t^{\a+1}}{\bar\r^2}\phi_y^2dyd\tau
+C (\f{\e}{\nu^{\g}\d})^2.
\end{array}\end{equation}
Recalling \eqref{assump2} and Cauchy inequality, it holds that
\begin{equation}\label{J8}\begin{array}{ll}
J_8&\dy=C|\int_0^\tau\int_{\mathbf{R}}\mu(\bar\t)\f{\phi_y}{\r^2}\big(\mu(\t)-\mu(\bar\t) \big)\psi_{yy} dyd\tau|\\[3mm]
&\dy\le \f{1}{16} \int_0^\tau\int_{\mathbf{R}}\f{\bar\t^{\a+1}}{\bar\r^2}\phi_y^2dyd\tau
+C \int_0^\tau\int_{\mathbf{R}}\f{\bar\t^{\a}}{\bar\r}\psi_{yy}^2\z^2~\bar\t^{2\a-3}\bar\r^{-1} dyd\tau\\[3mm]
&\dy\le \f{1}{16} \int_0^\tau\int_{\mathbf{R}}\f{\bar\t^{\a+1}}{\bar\r^2}\phi_y^2dyd\tau
+C\nu^{2-3\g}\sup_{[0,\tau_1(\e)]}\|\z\|\sup_{[0,\tau_1(\e)]}\|\z_y\|
\int_0^\tau\int_{\mathbf{R}}\f{\bar\t^{\a}}{\bar\r}\psi_{yy}^2 dyd\tau\\[3mm]
&\dy\le \f{1}{16} \int_0^\tau\int_{\mathbf{R}}\f{\bar\t^{\a+1}}{\bar\r^2}\phi_y^2dyd\tau
+C\nu^{1-3\g}\sup_{[0,\tau_1(\e)]}\|\sqrt{\bar\r^{2-\g}}\z\|\sup_{[0,\tau_1(\e)]}\|\z_y\|
\int_0^\tau\int_{\mathbf{R}}\f{\bar\t^{\a}}{\bar\r}\psi_{yy}^2 dyd\tau\\[3mm]
&\dy\le \f{1}{16} \int_0^\tau\int_{\mathbf{R}}\f{\bar\t^{\a+1}}{\bar\r^2}\phi_y^2dyd\tau
+C\nu^{1-3\g}\e^{\f16}
\int_0^\tau\int_{\mathbf{R}}\f{\bar\t^{\a}}{\bar\r}\psi_{yy}^2 dyd\tau.
\end{array}\end{equation}
Similarly,
\begin{equation}\label{J9}\begin{array}{ll}
J_9\dy=C|\a\int_0^\tau\int_{\mathbf{R}}\t^{\a-1}\z_y\psi_y\mu(\bar\t)\f{\phi_y}{\r^2} ~dyd\tau|\\[3mm]
\dy \le C\nu^{-\f{\g}{2}}\int_0^t\|\sqrt{\f{\bar\t^{2\a}}{\bar\r^3}}\phi_y\|
\|\sqrt{\bar\t^{\a-1}}\z_y \| \|\psi_y\|^{1/2} \|\psi_{yy}\|^{1/2} d\tau\\[3mm]
\dy \le C\nu^{-\f{2\g+\a(\g-1)}{4}}\int_0^t\|\sqrt{\f{\bar\t^{2\a}}{\bar\r^3}}\phi_y\|
\|\sqrt{\bar\t^{\a-1}}\z_y \| \|\psi_y\|^{1/2} \|\sqrt{\f{\bar\t^{\a}}{\bar\r}}\psi_{yy}\|^{1/2} d\tau\\[3mm]
\dy \le \nu^{2\a(\g-1)}|\ln\e|^{-1} \int_0^\tau\int_{\mathbf{R}}\f{\bar\t^{\a}}{\bar\r}\psi_{yy}^2 dyd\tau\\[3mm]
\dy\quad+C\nu^{-\f{2\g+3\a(\g-1)}{3}}|\ln\e|^{\f13}\int_0^\tau\|\sqrt{\f{\bar\t^{2\a}}{\bar\r^3}}\phi_y\|^{4/3}
\|\sqrt{\bar\t^{\a-1}}\z_y \|^{4/3} \|\psi_y\|^{2/3} d\tau\\[3mm]
\dy \le \nu^{2\a(\g-1)}|\ln\e|^{-1} \int_0^\tau\int_{\mathbf{R}}\f{\bar\t^{\a}}{\bar\r}\psi_{yy}^2 dyd\tau
+\f{1}{16}\sup_{[0,\tau_1(\e)]}\|\sqrt{\f{\bar\t^{2\a}}{\bar\r^3}}\phi_y \|^2\\[3mm]
\dy\quad+C \nu^{-2\g-3\a(\g-1)}|\ln\e|\Big(\int_0^\tau\|\sqrt{\bar\t^{\a-1}}\z_y \|^{4/3} \|\psi_y\|^{2/3} d\tau \Big)^3\\[3mm]
\dy \le \nu^{2\a(\g-1)}|\ln\e|^{-1} \int_0^\tau\int_{\mathbf{R}}\f{\bar\t^{\a}}{\bar\r}\psi_{yy}^2 dyd\tau
+\f{1}{16}\sup_{[0,\tau_1(\e)]}\|\sqrt{\f{\bar\t^{2\a}}{\bar\r^3}}\phi_y \|^2\\[3mm]
\dy\quad+C \nu^{-2\g-4\a(\g-1)}|\ln\e|\Big(\int_0^\tau\|\sqrt{\bar\t^{\a-1}}\z_y \|^2+ \|\sqrt{\bar\t^{\a}}\psi_y\|^2 d\tau \Big)^3\\[3mm]
\dy \le \nu^{2\a(\g-1)}|\ln\e|^{-1} \int_0^\tau\int_{\mathbf{R}}\f{\bar\t^{\a}}{\bar\r}\psi_{yy}^2 dyd\tau
+\f{1}{16}\sup_{[0,\tau_1(\e)]}\|\sqrt{\f{\bar\t^{2\a}}{\bar\r^3}}\phi_y \|^2
+C \nu^{-2\g-4\a(\g-1)}|\ln\e|\e\\[3mm]
\dy \le \nu^{2\a(\g-1)}|\ln\e|^{-1} \int_0^\tau\int_{\mathbf{R}}\f{\bar\t^{\a}}{\bar\r}\psi_{yy}^2 dyd\tau
+\f{1}{16}\sup_{[0,\tau_1(\e)]}\|\sqrt{\f{\bar\t^{2\a}}{\bar\r^3}}\phi_y \|^2+\e^{\f13},
\end{array}\end{equation}
where we have used the fact that
$$
C \nu^{-2\g-4\a(\g-1)}|\ln\e|\e=C\e^{1-(2\g+4\a(\g-1))a}|\ln\e|^{1-2\g-4\a(\g-1)}\le \e^{\f13},~~~{\rm if}~~\e\ll1.
$$
Combining \eqref{J1}-\eqref{J9} yields that
\begin{equation}\label{step21}\begin{array}{ll}
&\dy\sup_{\tau\in[0,\tau_1(\e)]}\int_{\mathbf{R}}
  \Big(\frac{\bar\t^{2\a}}{\bar\r^3}\phi_y^2+\bar\r^{\g-2}\phi^2+ \bar\r\psi^2+
  \bar\r^{2-\g}\zeta^2\Big)(\tau,y)dy\\[5mm]
+&\dy \int_0^{\tau_1(\e)} \int_{\mathbf{R}} \Big[\bar u_y\Big(\bar\r^{\g-2}\phi^2+ \bar\r\psi^2+
\bar\r^{2-\g}\zeta^2\Big) +\bar \theta^{\alpha}\psi_y^2+
\bar \theta^{\alpha-1}\zeta_y^2+\frac{\bar\theta^{\alpha+1}}{\bar\r^2}\phi_y^2\Big]dyd\tau\\[4mm]
\leq &\dy C\big(\nu^{2\a(\g-1)}|\ln\e|^{-1}+\nu^{1-3\g}\e^{\f16}\big)\int_0^{\tau_1(\e)}\int_{\mathbf{R}}
\frac{\bar\theta^{\alpha}}{\bar\r}\psi_{yy}^2 dyd\tau+C\e^{\f13}.
\end{array}\end{equation}
In particular, it holds that
\begin{equation}\label{step22}\begin{array}{ll}
&\dy\sup_{\tau\in[0,\tau_1(\e)]}\int_{\mathbf{R}}
  \frac{\bar\t^{2\a}}{\bar\r^3}\phi_y^2dy+
  \int_0^{\tau_1(\e)}\int_{\mathbf{R}}\frac{\bar\theta^{\alpha+1}}{\bar\r^2}\phi_y^2dyd\tau\\[4mm]
\leq &\dy C\big(\nu^{2\a(\g-1)}|\ln\e|^{-1}+\nu^{1-3\g}\e^{\f16}\big)\int_0^{\tau_1(\e)}\int_{\mathbf{R}}
\frac{\bar\theta^{\alpha}}{\bar\r}\psi_{yy}^2 dyd\tau+C\e^{\f13}.
\end{array}\end{equation}

\underline{\it Step 3. }\quad In the following, we estimate
$\dy\sup_\tau\|\psi_y\|$. For this, rewrite $\eqref{mass}_2$ as
\begin{equation}\label{remome}
\r\psi_\tau+\r u\psi_y+(\g-1)\big(\t\phi_y+\r\z_y\big)-\mu(\t)\psi_{yy}=-\bar g+\mu(\t)_y u_y,
\end{equation}
where
\begin{equation}\label{gbar}
\bar g=-\mu(\t)\bar u_{yy}+\r\psi \bar u_y+(\g-1)\big(\bar\r_y\z-\f{\bar\r_y\bar\t\phi}{\bar\r}\big).
\end{equation}
Multiplying \eqref{remome} by $-\psi_{yy}/\r$ gives
\begin{equation}
\begin{array}{l}
\dy(\f{\psi^2_y}{2})_\tau-(\psi_\tau\psi_y+u \f{\psi^2_y}{2})_y+\f12\bar u_y\psi^2_y+\f{\mu(\t)}{\r}\psi^2_{yy}\\
\dy=-\f{\psi^3_y}{2}+(\g-1)\big(\t\phi_y+\r\z_y\big)\f{\psi_{yy}}{\r}
+\bar g\f{\psi_{yy}}{\r}-\mu(\t)_y u_y\f{\psi_{yy}}{\r}.
\end{array}
\end{equation}
Integrating the above equation over ${\mathbf{R}}^1\times[0,\tau]$ yields
\begin{equation}\label{intmv}\begin{array}{ll}
&\dy\int_{\mathbf{R}}  \frac{\psi^2_y}{2}dy
+\int^\tau_{0}\int_{\mathbf{R}}\Big (\frac{{\bar u}_y\psi_y^2}{2}
  +\frac{\mu(\t)}{\r}\psi_{yy}^2\Big)dyd\tau\\[3mm]
=&\dy\int^\tau_{0}\int_{\mathbf{R}}\Big\{
\bar g\frac{\psi_{yy}}{\r}-\f{\psi_y^3}{2}+(\g-1)\big(\t\phi_y+\r\z_y \big)\f{\psi_{yy}}{\r}
-\mu(\t)_y u_y\f{\psi_{yy}}{\r} \Big\}dyd\tau.
 \end{array}\end{equation}
First, it follows from \eqref{g2} that
\begin{equation}\label{3.1}
\begin{array}{ll}
&\dy\Big|\int^\tau_{0}\int_{\mathbf{R}}
\bar g\frac{\psi_{yy}}{\r}dyd\tau\Big|\\[3mm]
\leq &\dy
\frac{1}{8}\int^\tau_{0}\int_{\mathbf{R}}
\f{\bar\t^{\a}}{\bar\r}\psi^2_{yy}dyd\tau+
C\Big|\int^\tau_{0}\int_{\mathbf{R}}\bar\t^{-\a}\bar\r^{-1}\bar g^2 dyd\tau\Big|\\[3mm]
\leq&\dy \frac{1}{8}\int^\tau_{0}\int_{\mathbf{R}}
\f{\bar\t^{\a}}{\bar\r}\psi^2_{yy}dyd\tau
+ C\nu^{-2\a(\g-1)}\int_0^\tau\int_{\mathbf{R}}\bar\t^{\a-1}\bar\r^{-2}g^2 dyd\tau\\[3mm]
\leq&\dy \frac{1}{8}\int^\tau_{0}\int_{\mathbf{R}}
\f{\bar\t^{\a}}{\bar\r}\psi^2_{yy}dyd\tau
+ C\nu^{-2\a(\g-1)-\g}\f{\e}{\d}\e^{\f13}
\le \frac{1}{8}\int^\tau_{0}\int_{\mathbf{R}}
\f{\bar\t^{\a}}{\bar\r}\psi^2_{yy}dyd\tau+\e^{\f13},
\end{array}
\end{equation}
where in the last inequality  we used the fact that
$$
C\nu^{-2\a(\g-1)-\g}\f{\e}{\d}=C\e^{1-a(\g+2\a(\g-1)+1)}|\ln\e|^{-2\a(\g-1)-\g}
\le C\e^{\f12}|\ln\e|^{-2\a(\g-1)-\g}\le 1,\quad {\rm
if}~~ \e\ll1.
$$
Furthermore, we can compute that
\begin{equation}\label{3.2}
\begin{array}{ll}
\dy\Big|\int^\tau_{0}\int_{\mathbf{R}}
\frac{\psi_y^3}{2}dyd\tau\Big| \leq
C\int^\tau_{0}\|\psi_{yy}\|^{\f12}\|\psi_y\|^{\f52} d\tau
\le C\nu^{-\a(\g-1)/4}\int^\tau_{0}\|\sqrt{\f{\bar\t^{\a}}{\bar\r}}\psi_{yy}\|^{\f12}\|\psi_y\|^{\f52}\\[3mm]
\qquad \dy \leq
\frac{1}{8}\int^\tau_{0}\int_{\mathbf{R}}\frac{\bar\t^{\a}}{\bar\r}\psi_{yy}^{2}
dyd\tau+C\nu^{-\a(\g-1)/3}\int^\tau_{0} \|\psi_y\|^{\f{10}{3}}d\tau\\[3mm]
\qquad \dy \leq
\frac{1}{8}\int^\tau_{0}\int_{\mathbf{R}}\frac{\bar\t^{\a}}{\bar\r}\psi_{yy}^{2}
dyd\tau+C\nu^{-4\a(\g-1)/3}\sup_{\tau\in[0,\tau_1]}\|\psi_y\|^{\f43}\int^\tau_{0}\int_{\mathbf{R}}
\bar\t^{\a}\psi_y^2 dyd\tau\\[3mm]
\qquad \dy \leq
\frac{1}{8}\int^\tau_{0}\int_{\mathbf{R}}\frac{\bar\t^{\a}}{\bar\r}\psi_{yy}^{2}
dyd\tau+C\nu^{-4\a(\g-1)/3}\e^{\f13},
\end{array}\end{equation}
where in the last inequality we have used the a priori assumptions \eqref{assump2}.
By Cauchy inequality, one has
\begin{equation}\label{3.3}
\begin{array}{ll}
&\dy\Big|\int^\tau_{0}\int_{\mathbf{R}}
(\t\phi_y+\r\z_y)\f{\psi_{yy}}{\r}dyd\tau\Big|\\[3mm]
&\dy \leq
\frac{1}{8}\int^\tau_{0}\int_{\mathbf{R}}\frac{\bar\t^{\a}}{\bar\r}\psi_{yy}^{2}
+C\int^\tau_{0}\int_{\mathbf{R}}
\big(\f{\bar\t^{\a+1}}{\bar\r^2}\phi_y^2+\bar\t^{\a-1}\z^2_y \big)\bar\r^{\g-2\a(\g-1)}
dyd\tau\\[3mm]
&\dy \le\frac{1}{8}\int^\tau_{0}\int_{\mathbf{R}}\frac{\bar\t^{\a}}{\bar\r}\psi_{yy}^{2}
+C\nu^{-2\a(\g-1)} \int^\tau_{0}\int_{\mathbf{R}}
\f{\bar\t^{\a+1}}{\bar\r^2}\phi_y^2 dyd\tau+C\nu^{-2\a(\g-1)}\e^{\f13}.
\end{array}\end{equation}
Finally, one has
\begin{equation}\label{3.4}
\begin{array}{ll}
&\dy\Big|\int^\tau_{0}\int_{\mathbf{R}}
\mu(\t)_y u_y\f{\psi_{yy}}{\r}dyd\tau\Big|\\[3mm]
\leq &\dy  C \int^\tau_{0}\int_{\mathbf{R}}
\bar\t^{\a-1}\bar\r^{-1}|\psi_{yy}|\big(|\psi_y\z_y|+\bar u_y(\bar\t^{1/2}|\psi_y|+|\z_y|)+\bar\t^{1/2}\bar u^2_y \big)
dyd\tau:=\sum_{i=1}^3K_i.
\end{array}\end{equation}
The terms $K_i~(i=1,2,3)$ will be estimated as follows,
\begin{equation}\label{K1}\begin{array}{ll}
K_1&\dy=C\int_0^\tau\int_{\mathbf{R}}
\bar\t^{\a-1}\bar\r^{-1}|\psi_{yy}\psi_y\z_y| dyd\tau \le C\nu^{-\g/2}\int_0^\tau\|\sqrt{\f{\bar\t^{\a}}{\bar\r}}\psi_{yy}\| \|\sqrt{\bar\t^{\a-1}}\z_y\| \|\psi_y\|^{\f12} \|\psi_{yy}\|^{\f12}d\tau\\[3mm]
&\dy\le C\nu^{-\f{2\g+\a(\g-1)}{4}}\int_0^\tau\|\sqrt{\f{\bar\t^{\a}}{\bar\r}}\psi_{yy}\|^{\f32} \|\sqrt{\bar\t^{\a-1}}\z_y\| \|\psi_y\|^{\f12} d\tau\\[3mm]
&\dy \le\frac{1}{8}\int^\tau_{0}\int_{\mathbf{R}}\frac{\bar\t^{\a}}{\bar\r}\psi_{yy}^{2}
+C\nu^{-2\g-\a(\g-1)}\int_0^\tau \|\sqrt{\bar\t^{\a-1}}\z_y\|^4 \|\psi_y\|^2 d\tau  \\[3mm]
&\dy \le\frac{1}{8}\int^\tau_{0}\int_{\mathbf{R}}\frac{\bar\t^{\a}}{\bar\r}\psi_{yy}^{2} dyd\tau
+C\nu^{1-3\g-\a(\g-1)}\sup_{[0,\tau_1(\e)]}\|\z_y\|^2 \sup_{[0,\tau_1(\e)]}\|\psi_y\|^2
\int^\tau_{0}\int_{\mathbf{R}}\bar\t^{\a-1}\z^2_y dyd\tau\\[3mm]
&\dy \le\frac{1}{8}\int^\tau_{0}\int_{\mathbf{R}}\frac{\bar\t^{\a}}{\bar\r}\psi_{yy}^{2} dyd\tau
+C\nu^{1-3\g-\a(\g-1)}\e^{\f13}\\[3mm]
&\dy \le\frac{1}{8}\int^\tau_{0}\int_{\mathbf{R}}\frac{\bar\t^{\a}}{\bar\r}\psi_{yy}^{2} dyd\tau
+\nu^{1-3\g-2\a(\g-1)}\e^{\f13}
\le\frac{1}{8}\int^\tau_{0}\int_{\mathbf{R}}\frac{\bar\t^{\a}}{\bar\r}\psi_{yy}^{2} dyd\tau+\e^{\f16+a},
\end{array}\end{equation}
where we used the fact that
$$
\nu^{1-3\g-2\a(\g-1)}\e^{\f16}=\e^a|\ln\e|^{1-3\g-2\a(\g-1)}\le\e^a,~~~{\rm if}~~\e\ll1.
$$
By Cauchy inequality, it holds that
\begin{equation}\label{K2}\begin{array}{ll}
K_2&\dy=C\int_0^\tau\int_{\mathbf{R}}
\bar\t^{\a-1}\bar\r^{-1}|\psi_{yy}| \bar u_y (\bar\t^{\f12}|\psi_y|+|\z_y|) dyd\tau\\[3mm]
&\dy \le\frac{1}{8}\int^\tau_{0}\int_{\mathbf{R}}\frac{\bar\t^{\a}}{\bar\r}\psi_{yy}^{2} dyd\tau
+C\int_0^\tau\int_{\mathbf{R}}
(\bar\t^{\a}\psi_y^2+\bar\t^{\a-1}\z_y^2)(\bar\t\bar\r)^{-1}\bar u_y^2
dyd\tau\\[3mm]
&\dy \le\frac{1}{8}\int^\tau_{0}\int_{\mathbf{R}}\frac{\bar\t^{\a}}{\bar\r}\psi_{yy}^{2} dyd\tau
+C\nu^{-\g}(\f{\e}{\d})^2\e^{\f13}
\le \frac{1}{8}\int^\tau_{0}\int_{\mathbf{R}}\frac{\bar\t^{\a}}{\bar\r}\psi_{yy}^{2} dyd\tau
+\e^{\f13}.
\end{array}\end{equation}
Recalling Lemma 2.3, one can get
\begin{equation}\label{K3}\begin{array}{ll}
K_3&\dy=C\int_0^\tau\int_{\mathbf{R}}
\bar\t^{\a-\f12}\bar\r^{-1}|\psi_{yy}|\bar u_y^2
dtd\tau\\[3mm]
&\dy \le\frac{1}{8}\int^\tau_{0}\int_{\mathbf{R}}\frac{\bar\t^{\a}}{\bar\r}\psi_{yy}^{2} dyd\tau
+C\int_0^\tau\int_{\mathbf{R}}
\bar\t^{\a-1}\bar\r^{-1}\bar u_y^4
dyd\tau\\[3mm]
&\dy \le\frac{1}{8}\int^\tau_{0}\int_{\mathbf{R}}\frac{\bar\t^{\a}}{\bar\r}\psi_{yy}^{2} dyd\tau
+C\nu^{-\g}\int_0^\tau \|\bar u_y\|^4_{L^4}
d\tau\\[3mm]
&\dy \le\frac{1}{8}\int^\tau_{0}\int_{\mathbf{R}}\frac{\bar\t^{\a}}{\bar\r}\psi_{yy}^{2} dyd\tau
+C\nu^{-\g}(\f{\e}{\d})^2
\le \frac{1}{8}\int^\tau_{0}\int_{\mathbf{R}}\frac{\bar\t^{\a}}{\bar\r}\psi_{yy}^{2} dyd\tau
+\e^{\f13}.
\end{array}\end{equation}
Substituting \eqref{3.1}-\eqref{K3} into \eqref{intmv}, it holds that
\begin{equation}\label{psiy}
\begin{array}{ll}
\dy\sup_{\tau\in[0,\tau_1(\e)]}\int_{\mathbf{R}}\psi_y^2dy+ \int^{\tau_1(\e)}_{0}\int_{\mathbf{R}}\Big
({\bar u}_y\psi_y^2
  +\frac{\bar\t^{\a}}{\bar\r}\psi_{yy}^2\Big)dyd\tau\\
 \qquad\qquad \dy \leq\ C\nu^{-2\a(\g-1)}\int^{\tau_1(\e)}_{0}\int_{\mathbf{R}}
\f{\bar\t^{\a+1}}{\bar\r^2}\phi^2_y
dyd\tau+\e^{\f16+a}.
\end{array}\end{equation}
Then substituting \eqref{psiy} into \eqref{step22}, one can get
\begin{equation}\label{psiy1}
\begin{array}{ll}
\dy\sup_{\tau\in[0,\tau_1(\e)]}\int_{\mathbf{R}}
  \frac{\bar\t^{2\a}}{\bar\r^3}\phi_y^2dy+
  \int^{\tau_1(\e)}_{0}\int_{\mathbf{R}}\frac{\bar\theta^{\alpha+1}}{\bar\r^2}\phi_y^2dyd\tau\\[4mm]
  \leq \dy \big(\nu^{2\a(\g-1)}|\ln\e|^{-1}+C\nu^{1-3\g}\e^{\f16}\big)\int_0^{\tau_1(\e)}\int_{\mathbf{R}}
\frac{\bar\theta^{\alpha}}{\bar\r}\psi_{yy}^2 dyd\tau+C\e^{\f13}\\[4mm]
\le \dy C\big(|\ln\e|^{-1}+\e^a \big)\int^{\tau_1(\e)}_{0}\int_{\mathbf{R}}
\f{\bar\t^{\a+1}}{\bar\r^2}\phi^2_y
dyd\tau+\nu^{2\a(\g-1)}\e^{\f16+a}|\ln\e|^{-1}+C\nu^{1-3\g}\e^{\f13+a}+C\e^{\f13}.
\end{array}\end{equation}
Note that
$$
\nu^{2\a(\g-1)}\e^{\f16+a}|\ln\e|^{-1}=\e^{\f13-3a\g+a}|\ln\e|^{2\a(\g-1)-1} \leq  \e^{\f13-3a\g}|\ln\e|^{-3\g},\quad {\rm
if}~~ \e\ll1,
$$
and
$$
C\nu^{1-3\g}\e^{\f13+a}=C\e^{\f13-3a\g+2a}|\ln\e|^{1-3\g}\leq  \e^{\f13-3a\g}|\ln\e|^{-3\g},\quad {\rm
if}~~ \e\ll1,
$$
one can get from \eqref{psiy} and \eqref{psiy1} that
\begin{equation}\label{phiy}\begin{array}{ll}
&\dy\sup_{\tau\in[0,\tau_1(\e)]}\int_{\mathbf{R}}
  \frac{\bar\t^{2\a}}{\bar\r^3}\phi_y^2dy+
  \int^{\tau_1(\e)}_{0}\int_{\mathbf{R}}\frac{\bar\theta^{\alpha+1}}{\bar\r^2}\phi_y^2dyd\tau
\leq  \e^{\f13-3a\g}|\ln\e|^{-3\g},\quad {\rm
if}~~ \e\ll1.
\end{array}\end{equation}
Meanwhile, it holds that
\begin{equation}\label{psiy2}
\begin{array}{ll}
&\dy\sup_{\tau\in[0,\tau_1(\e)]}\int_{\mathbf{R}}\psi_y^2dy+ \int^{\tau_1(\e)}_{0}\int_{\mathbf{R}}\Big
({\bar u}_y\psi_y^2
  +\frac{\bar\t^{\a}}{\bar\r}\psi_{yy}^2\Big)dyd\tau\\[4mm]
&\dy\leq C\nu^{-2\a(\g-1)}\e^{\f13-3a\g}|\ln\e|^{-3\g}+\e^{\f16+a}
=C\e^{\f16}|\ln\e|^{-3\g-2\a(\g-1)}+\e^{\f16+a}\le\e^{\f16}.
\end{array}\end{equation}
\underline{\it Step 4. }\quad Finally, we estimate $\dy\sup_\tau\|\z_y\|$. For this, rewrite $\eqref{mass}_3$ as
\begin{equation}\label{reenergy}
\r\z_\tau+\r u\z_y+(\g-1)\r\t\psi_y-\k(\t)\z_{yy}=-\bar h+\k(\t)_y \t_y+\mu(\t)u^2_y,
\end{equation}
where
\begin{equation}\label{hbar}
\bar h=-\k(\t)\bar \t_{yy}+\r\psi \bar \t_y+(\g-1)\r\z\bar u_y.
\end{equation}
Multiplying \eqref{reenergy} by $-\z_{yy}/\r$ gives
\begin{equation}
\begin{array}{l}
\dy(\f{\z^2_y}{2})_\tau-(\z_\tau\z_y+u \f{\z^2_y}{2})_y+\f12\bar u_y\z^2_y+\f{\k(\t)}{\r}\z^2_{yy}\\
\dy=\bar h\f{\z_{yy}}{\r}+(\g-1)\t\psi_y\z_{yy}-\f{\psi_y\z^2_y}{2}
-\big(\k(\t)_y \t_y+\mu(\t)u^2_y \big)\f{\z_{yy}}{\r}.
\end{array}
\end{equation}
Integrating the above equation over ${\mathbf{R}}^1\times[0,\tau]$ yields
\begin{equation}\label{intzy}\begin{array}{ll}
&\dy\int_{\mathbf{R}}  \frac{\z^2_y}{2}dy
+\int^\tau_{0}\int_{\mathbf{R}}\Big (\frac{{\bar u}_y\z_y^2}{2}
  +\frac{\k(\t)}{\r}\z_{yy}^2\Big)dyd\tau\\[3mm]
=&\dy\int^\tau_{0}\int_{\mathbf{R}}\Big\{
\bar h\frac{\z_{yy}}{\r}+(\g-1)\t\psi_y\z_{yy}-\f{\psi_y\z^2_y}{2}
-\big(\k(\t)_y \t_y+\mu(\t)u^2_y \big)\f{\z_{yy}}{\r} \Big\}dyd\tau.
 \end{array}\end{equation}
First,
\begin{equation}\label{4.1}
\dy\int_0^\tau\int_{\mathbf{R}}\bar h\f{\z_{yy}}{\r}
dyd\tau\le \f{1}{8}\int_0^\tau\int_{\mathbf{R}}
\f{\bar\t^{\a}}{\bar\r}\z^2_{yy}
dyd\tau+C \int_0^\tau\int_{\mathbf{R}}\bar\t^{-\a}\bar\r^{-1}|\bar h|^2
dyd\tau.
\end{equation}
Recalling \eqref{rup},  \eqref{hbar},  and $(i)$ in Lemma 2.3, one can get
\begin{equation}\begin{array}{ll}
|\bar h|&\le\dy  C\big(\bar\t^{\a}|\bar \t_{yy}|+|\bar\r\bar \t_y\psi|+|\bar u_y\bar\r\z| \big)\\[2mm]
&\dy\le C\Big\{\bar\t^{\a}(\bar\t^{\f12}|\bar u_{yy}|+|\bar u_y|^2)
+\bar u_y(|\bar\t^{\f12}\bar\r\psi|+|\bar\r\z|) \Big\}.
\end{array}\end{equation}
So the last term in \eqref{4.1} can be estimated by
\begin{equation}\begin{array}{ll}
&\dy\int_0^\tau\int_{\mathbf{R}}\bar\t^{\a}\bar\r^{-1}|\bar h|^2 dyd\tau\\[2mm]
&\dy\le C\nu^{-1}\int_0^\tau(\|\bar u_{yy}\|^2+ \|\bar u_y\|^4_{L^4})d\tau
+C\int_0^\tau\int_{\mathbf{R}}
\bar u_y(\bar\r\psi^2+\bar\r^{2-\g}\z^2)\bar\t^{1-\a}\bar u_y
dyd\tau\\[3mm]
&\dy\le C\nu^{-1}(\f{\e}{\d})^2+C\nu^{-\a(\g-1)}\f{\e}{\d}\e^{\f13} \le \e^{\f13}.
\end{array}\label{h2}\end{equation}
By Cauchy inequality, it holds that
\begin{equation}\label{4.21}
\begin{array}{ll}
\dy\Big|\int^\tau_{0}\int_{\mathbf{R}}
\t\psi_y\z_{yy}dyd\tau\Big|
\dy \leq
\frac{1}{8}\int^\tau_{0}\int_{\mathbf{R}}\frac{\bar\t^{\a}}{\bar\r}\z_{yy}^{2}
+C\int^\tau_{0}\int_{\mathbf{R}}
\bar\t^{\a}\psi^2_y\bar\t^{2-2\a}\bar\r
dyd\tau\\[3mm]
\qquad\qquad\dy \le\frac{1}{8}\int^\tau_{0}\int_{\mathbf{R}}\frac{\bar\t^{\a}}{\bar\r}\z_{yy}^{2}
+C\nu^{-2\a(\g-1)}\e^{\f13}.
\end{array}\end{equation}
By Sobolev inequality and \eqref{psiy}, it holds that
\begin{equation}\label{4.2}
\begin{array}{ll}
&\dy\Big|\int^\tau_{0}\int_{\mathbf{R}}
\frac{\psi_y\z^2_y}{2}dyd\tau\Big|\leq
C\int^\tau_{0}\|\psi_{y}\|\|\z_y\|^2_{L^4} d\tau
\le C\int^\tau_{0}\|\psi_{y}\|\|\z_y\|^{\f32}\|\z_{yy}\|^{\f12} d\tau\\[3mm]
&\qquad\dy \leq
\frac{1}{8}\int^\tau_{0}\int_{\mathbf{R}}\frac{\bar\t^{\a}}{\bar\r}\z_{yy}^{2}
dyd\tau+C\nu^{-\a(\g-1)/3}\int^\tau_{0} \|\psi_y\|^{\f{4}{3}}\|\z_y\|^2d\tau\\[3mm]
&\qquad\dy \leq
\frac{1}{8}\int^\tau_{0}\int_{\mathbf{R}}\frac{\bar\t^{\a}}{\bar\r}\z_{yy}^{2}
dyd\tau+C\nu^{-4\a(\g-1)/3}\sup_{\tau\in[0,\tau_1]}\|\psi_y\|^{\f43}\int^\tau_{0}\int_{\mathbf{R}}
\bar\t^{\a-1}\z_y^2 dyd\tau\\[4mm]
&\qquad\dy \leq
\frac{1}{8}\int^\tau_{0}\int_{\mathbf{R}}\frac{\bar\t^{\a}}{\bar\r}\z_{yy}^{2}
dyd\tau+C\nu^{-4\a(\g-1)/3}\e^{\f13} \leq
\frac{1}{8}\int^\tau_{0}\int_{\mathbf{R}}\frac{\bar\t^{\a}}{\bar\r}\psi_{yy}^{2}
dyd\tau+\e^{\f19},
\end{array}\end{equation}
where we have used the a priori assumptions \eqref{assump2}, and in the last inequality we have used the fact that
$$
C\nu^{-4\a(\g-1)/3}\e^{\f13}=C\e^{\f{1-4\a(\g-1)a}{3}}|\ln\e|^{-\f{4\a(\g-1)}{3}}\le \e^{\f19},~~{\rm if} ~~\e\ll1.
$$
By Cauchy inequality, we have
\begin{equation}\label{4.3}
\begin{array}{ll}
&\dy\Big|\int^\tau_{0}\int_{\mathbf{R}}
(\k(\t)_y\t_y+\mu(\t)u^2_y)\f{\z_{yy}}{\r}dyd\tau\Big|\\[3mm]
\le &\dy C\int^\tau_{0}\int_{\mathbf{R}}
\bar\t^{\a}\bar\r^{-1}|\z_{yy}|(\bar\t^{-1}\z^2_y+\psi^2_y+\bar u^2_y)
dyd\tau:=\sum^3_{i=1}L_i.
\end{array}\end{equation}
Now we estimate the terms on the right-hand side of \eqref{4.3} one by one. By Sobolev inequality, it holds that
\begin{equation}\label{L1}\begin{array}{ll}
L_1&\dy=\int_0^\tau\int_{\mathbf{R}}
\bar\t^{\a-1}\bar\r^{-1}|\z_{yy}\z^2_y| dyd\tau
 \le C\nu^{1/2-\g}\int_0^\tau\|\sqrt{\f{\bar\t^{\a}}{\bar\r}}\z_{yy}\| \|\z_y\|^2_{L^4}
d\tau\\[3mm]
&\dy \le C\nu^{1/2-\g}\int_0^\tau\|\sqrt{\f{\bar\t^{\a}}{\bar\r}}\z_{yy}\| \|\z_{yy}\|^{\f12} \|\z_y\|^{\f32}
d\tau
 \le C\nu^{\f{2-4\g-\a(\g-1)}{4}}\int_0^\tau\|\sqrt{\f{\bar\t^{\a}}{\bar\r}}\z_{yy}\|^{\f32} \|\z_y\|^{\f32}
d\tau\\[3mm]
&\dy \le \f18\int^\tau_{0}\int_{\mathbf{R}}\frac{\bar\t^{\a}}{\bar\r}\z_{yy}^{2}
dyd\tau+C\nu^{2-4\g-\a(\g-1)}\int_0^\tau\|\z_y\|^6 d\tau\\[3mm]
&\dy \le \f18\int^\tau_{0}\int_{\mathbf{R}}\frac{\bar\t^{\a}}{\bar\r}\z_{yy}^{2}
dyd\tau+C\nu^{2-4\g-2\a(\g-1)}\sup_{\tau\in[0,\tau_1(\e)]}\|\z_y\|^4
\int_0^\tau\int_{\mathbf{R}}\bar\t^{\a-1}\z_y^2 dyd\tau\\[4mm]
&\dy \le \f18\int^\tau_{0}\int_{\mathbf{R}}\frac{\bar\t^{\a}}{\bar\r}\z_{yy}^{2}
dyd\tau+C\nu^{2-4\g-2\a(\g-1)}\e^{\f13}\le \f18\int^\tau_{0}\int_{\mathbf{R}}\frac{\bar\t^{\a}}{\bar\r}\z_{yy}^{2}
dyd\tau+\e^{\f19},
\end{array}\end{equation}
where in the last inequality we have used the fact that
$$
C\nu^{2-4\g-2\a(\g-1)}\e^{\f13}=C\e^{\f13+a(2-4\g-2\a(\g-1))}|\ln\e|^{2-4\g-2\a(\g-1)}\le \e^{\f19},~~{\rm if}~~\e\ll1.
$$
Similarly, one has
\begin{equation}\label{L2}\begin{array}{ll}
L_2&\dy=\int_0^\tau\int_{\mathbf{R}}
\bar\t^{\a}\bar\r^{-1}|\z_{yy}\psi^2_y| dyd\tau\le C\nu^{-\f12}\int_0^\tau\|\sqrt{\frac{\bar\t^{\a}}{\bar\r}}\z_{yy}\| \|\psi_y\|^2_{L^4}
d\tau
\\[3mm]
&\dy \le C\nu^{-\f12}\int_0^\tau\|\sqrt{\frac{\bar\t^{\a}}{\bar\r}}\z_{yy}\| \|\psi_{yy}\|^{\f12} \|\psi_y\|^{\f32}
d\tau \le C\nu^{-\f{2+\a(\g-1)}{4}}\int_0^\tau\|\sqrt{\f{\bar\t^{\a}}{\bar\r}}\z_{yy}\|
\|\sqrt{\f{\bar\t^{\a}}{\bar\r}}\psi_{yy}\|^{\f12} \|\psi_y\|^{\f32}
d\tau\\[3mm]
&\dy \le \f18\int^\tau_{0}\int_{\mathbf{R}}\frac{\bar\t^{\a}}{\bar\r}\z_{yy}^{2}
dyd\tau+C\nu^{-\f{2+\a(\g-1)}{2}}\int_0^\tau\|\sqrt{\f{\bar\t^{\a}}{\bar\r}}\psi_{yy}\| \|\psi_y\|^3
d\tau\\[3mm]
&\dy \le \f18\int^\tau_{0}\int_{\mathbf{R}}(\frac{\bar\t^{\a}}{\bar\r}\z_{yy}^{2}+\frac{\bar\t^{\a}}{\bar\r}\psi_{yy}^{2})
dyd\tau+C\nu^{-2-\a(\g-1)}\int_0^\tau \|\psi_y\|^6
d\tau\\[3mm]
&\dy \le \f18\int^\tau_{0}\int_{\mathbf{R}}\frac{\bar\t^{\a}}{\bar\r}\z_{yy}^{2}
dyd\tau +\f18\e^{\f19} +C\nu^{-2-2\a(\g-1)}\sup_{[0,\tau_1(\e)]}\|\psi_y\|^4
\int_0^\tau\int_{\mathbf{R}}\bar\t^{\a}\psi_y dyd\tau\\[4mm]
&\dy \le \f18\int^\tau_{0}\int_{\mathbf{R}}\frac{\bar\t^{\a}}{\bar\r}\z_{yy}^{2}
dyd\tau +\f18\e^{\f19} +C\nu^{-2-2\a(\g-1)}\e^{\f59} \le \f18\int^\tau_{0}\int_{\mathbf{R}}\frac{\bar\t^{\a}}{\bar\r}\z_{yy}^{2}
dyd\tau+\e^{\f19}.
\end{array}\end{equation}
Recalling Lemma 2.3 and using Cauchy inequality, it holds that
\begin{equation}\label{L3}\begin{array}{ll}
L_3&\dy=\int_0^\tau\int_{\mathbf{R}}
\bar\t^{\a}\bar\r^{-1}|\z_{yy}|\bar u_y^2
dtd\tau\le\frac{1}{8}\int^\tau_{0}\int_{\mathbf{R}}\frac{\bar\t^{\a}}{\bar\r}\z_{yy}^{2} dyd\tau
+C\int_0^\tau\int_{\mathbf{R}}
\bar\t^{\a}\bar\r^{-1}\bar u_y^4
dyd\tau\\[4mm]
&\dy \le\frac{1}{8}\int^\tau_{0}\int_{\mathbf{R}}\frac{\bar\t^{\a}}{\bar\r}\z_{yy}^{2} dyd\tau
+C\nu^{-1}\int_0^\tau \|\bar u_y\|^4_{L^4}
d\tau\le\frac{1}{8}\int^\tau_{0}\int_{\mathbf{R}}\frac{\bar\t^{\a}}{\bar\r}\z_{yy}^{2} dyd\tau
+\e^{\f19}.
\end{array}\end{equation}
Substituting \eqref{4.1}, \eqref{4.2}-\eqref{L3} into \eqref{intzy}, it holds that
\begin{equation}\label{zy}
\begin{array}{ll}
&\dy\sup_{\tau\in[0,\tau_1(\e)]}\int_{\mathbf{R}}\z_y^2dy+ \int^{\tau_1(\e)}_{0}\int_{\mathbf{R}}\Big
({\bar u}_y\psi_y^2
  +\frac{\bar\t^{\a}}{\bar\r}\z_{yy}^2\Big)dyd\tau \leq \e^{\f19}.
\end{array}\end{equation}
Therefore, \eqref{main1}, \eqref{main2} and \eqref{main3} can be derived directly from $(\ref{step1})$,
\eqref{phiy}, \eqref{psiy2} and \eqref{zy}.
It follows from \eqref{main1}-\eqref{main3} that if $\e$ is suitably small, then
\begin{equation*}
\begin{array}{ll}
\dy\sup_{0\leq\tau\leq\tau_1(\e)}\|\phi(\cdot,\tau)\|_{L^\infty} &\di
\leq\sqrt2\sup_{0\leq\tau\leq\tau_1(\e)}\|\phi(\cdot,\tau)\|^{1/2}\|\phi_y(\cdot,\tau)\|^{1/2}\\[3mm]
&\di \leq
C\sup_{0\leq\tau\leq\tau_1(\e)}
\Big(\nu^{-\g-2\a(\g-1)}\int_{\mathbf{R}}\bar\r^{\g-2}\phi^2dy
\int_{\mathbf{R}}\f{\bar\t^{2\a}}{\bar\r^3}\phi_y^2dy \Big)^{\f14}
\\[5mm]
&\di \leq C\Big(\nu^{-\g-2\a(\g-1)}\e^{\f13}\cdot\e^{\f13-3a\g}|\ln\e|^{-3\g}\Big)^{\f14}\\[5mm]
&\di =C\e^{\f{4\g+3\a(\g-1)}{2(18\g+12\a(\g-1))}}|\ln\e|^{-\f{2\g+\a(\g-1)}{2}}
 \le \e^{\f{1}{18\g+12\a(\g-1)}}=\e^a,
\end{array}
\end{equation*}
\begin{equation*}
\begin{array}{ll}
\dy\sup_{0\leq\tau\leq\tau_1(\e)}\|\psi(\cdot,\tau)\|_{L^\infty} &\di
\leq\sqrt2\sup_{0\leq\tau\leq\tau_1(\e)}\|\psi(\cdot,\tau)\|^{1/2}\|\psi_y(\cdot,\tau)\|^{1/2}\\[3mm]
&\di \leq
C\sup_{0\leq\tau\leq\tau_1(\e)}
\Big(\nu^{-1}\int_{\mathbf{R}}\bar\r\psi^2dy
\int_{\mathbf{R}}\psi_y^2dy \Big)^{\f14}
\\[5mm]
&\di \leq C\Big(\nu^{-1}\e^{\f13}\cdot \e^{\f19}\Big)^{\f14}
=C\e^{\f19-\f{a}{4}}|\ln\e|^{-\f14}\\[3mm]
&\di=C\e^{\f{24\g+16\a(\g-1)-3}{12(18\g+12\a(\g-1))}}|\ln\e|^{-\f14}
\le \e^{\f{1}{18\g+12\a(\g-1)}}= \e^a,
\end{array}
\end{equation*}
and
\begin{equation*}
\begin{array}{ll}
\dy\sup_{0\leq\tau\leq\tau_1(\e)}\|\z(\cdot,\tau)\|_{L^\infty} &\di
\leq\sqrt2\sup_{0\leq\tau\leq\tau_1(\e)}\|\z(\cdot,\tau)\|^{1/2}\|\z_y(\cdot,\tau)\|^{1/2}\\[3mm]
&\di \leq
C\sup_{0\leq\tau\leq\tau_1(\e)}
\Big(\nu^{-2}\int_{\mathbf{R}}\bar\r^{2-\g}\z^2dy
\int_{\mathbf{R}}\z_y^2dy \Big)^{\f14}
\\[5mm]
&\di \leq C\Big(\nu^{-2}\e^{\f13}\cdot \e^{\f19}\Big)^{\f14}
=C\e^{\f19-\f{a}{2}}|\ln\e|^{-\f12}
\\[3mm]
&\di=C\e^{\f{12\g+8\a(\g-1)-3}{6(18\g+12\a(\g-1))}}|\ln\e|^{-\f12}
\le \e^{\f{\g-1}{18\g+12\a(\g-1)}}=\e^{(\g-1)a}.
\end{array}
\end{equation*}
Thus the a priori assumptions
\eqref{assump}-\eqref{assump2} are verified and the proof of Lemma \ref{len} is completed.
 \hfill $\Box$
\vspace{2mm}\\

It is note that the a priori estimates (\ref{main1})-(\ref{main3}) are better than the a priori assumptions (\ref{assump})-(\ref{assump2}) in the time interval $[0,\tau_1(\e)]$ with $\tau_1(\e)$ being the maximum existence time. Based on these a priori estimates, we can claim $\tau_1(\e)=\infty$. In fact, if $\tau_1(\e)<\infty$, then by again using the local existence at time $\tau=\tau_1(\e)$, we can find another time $\tau_2(\e)>\tau_1(\e)$ so that the solution satisfies the assumptions (\ref{assump})-(\ref{assump2}) in the time interval $[0,\tau_2(\e)]$ which contradicts the assumption that $\tau_1(\e)$ is the maximum existence time. Therefore we extend the local solution to the global one in $[0,\infty)$ for small but fixed $\e$.

\textbf{ Proof of Theorem ~\ref{thm1}:}\ \
 It remains to prove (\ref{cr}) with $a$ given in (\ref{a}).  From Lemma \ref{cut-off}, Lemma \ref{appu} (iii), \eqref{main1+}-\eqref{main3+}
  and recalling that
$\nu=\e^{a}|\ln\e|,~~\delta=\e^{a}$ in \eqref{mu}, it holds that for any given
positive constant $l$, there exists a constant $C_l>0$ which is
independent of $\epsilon$ such that

\begin{equation*}\begin{array}{ll}
&\dy\sup_{t\geq l}\| \r(\cdot,t)-\rho^{r_3}(\frac{\cdot}{t})\|_{L^\infty}\\
\leq&\dy\sup_{\tau\in[0,+\i)}\|
\phi(\cdot,\tau)\|_{L^\infty}+\sup_{t\geq
l}\|\bar\r(\cdot,t)-\rho^{r_3}_\nu(\frac{\cdot}{t})\|
_{L^\infty}+\sup_{t\geq l}\| \rho^{r_3}_\nu(\frac{\cdot}{t})- \rho^{r_3}(\frac{\cdot}{t})\|_{L^\infty}\\[4mm]
\leq &\dy C_l\big(\e^a+\d|\ln\d|+\nu \big)
\leq  C_l~\epsilon^{a}|\ln\epsilon |.
\end{array}\end{equation*}
About the convergence of the momentum, we have
\begin{equation*}\begin{array}{ll}
&\dy\sup_{t\geq l}\| m(\cdot,t)-m^{r_3}(\frac{\cdot}{t})\|_{L^\infty}\\
\leq&\dy C\sup_{\tau\in[0,+\i)}\Big(\|\psi(\cdot,\tau)\|_{L^\infty}+\|
\phi(\cdot,\tau)\|_{L^\infty}\Big)+\sup_{t\geq l}\Big(\|\bar
m(\cdot,t)-m^{r_3}_\nu(\frac{\cdot}{t})\|
_{L^\infty}+\| m^{r_3}_\nu(\frac{\cdot}{t})- m^{r_3}(\frac{\cdot}{t})\|_{L^\infty}\Big)\\
\leq &\dy C_l\big(\e^a+\delta|\ln\delta|+\nu\big)
\leq  C_l~\epsilon^{a} |\ln\epsilon |.
\end{array}\end{equation*}
About the convergence of the total energy, we can get
\begin{equation*}\begin{array}{ll}
&\dy\sup_{t\geq l}\| n(\cdot,t)-n^{r_3}(\frac{\cdot}{t})\|_{L^\infty}\\
\leq&\dy C\sup_{\tau\in[0,+\i)}\Big(\|\z(\cdot,\tau)\|_{L^\infty}+\|\phi(\cdot,\tau)\|_{L^\infty}\Big)+\sup_{t\geq
l}\Big(\|\bar\r\bar\t(\cdot,t)-n^{r_3}_\nu(\frac{\cdot}{t})\|
_{L^\infty}+\| n^{r_3}_\nu(\frac{\cdot}{t})- n^{r_3}(\frac{\cdot}{t})\|_{L^\infty} \Big)
\\[4mm]
\leq &\dy C_l\big(\e^a+\d|\ln\d|+\nu \big)
\leq  C_l~\epsilon^{a}|\ln\epsilon |,
\end{array}\end{equation*}
where we have used the fact that
$$
\di\sup_{\tau\in[0,+\i)}\|\z(\cdot,\tau)\|_{L^\infty}\le
C\e^{\f19-\f{a}{2}}|\ln\e|^{-\f12}\le \e^a.
$$
Then the proof of Theorem 1.1 is completed.
 \hfill $\Box$
\vspace{2mm}


\textbf{    }
\end{CJK}

\begin{thebibliography}{99}
\bibitem{B-B} S. Bianchini, A. Bressan, Vanishing viscosity
solutions of nonlinear hyperbolic systems,
 Ann. of Math. (2), 161 (2005), pp. 223-342.

\bibitem{Ch-C} S. Chapman, T.G. Cowling, The Mathematical Theory of Non-Uniform Gases, Cambridge University Press, 3rd edition, 1990.


\bibitem{chen-P} G. Q. Chen, M. Perepelitsa, Vanishing viscosity limit of the
Navier-Stokes equations to the Euler equations for compressible
fluid flow, Comm. Pure Appl. Math., 63, (2010), pp. 1469-1504.


\bibitem{good xin}
J. Goodman, Z. P. Xin, Viscous limits for piecewise smooth solutions
to systems of conservation laws, Arch. Ration. Mech. Anal., 121
(1992), pp. 235-265.


\bibitem{hoffliu}
D. Hoff, T. P. Liu, The inviscid limit for the Navier-Stokes
equations of compressible, isentropic flow with shock data, Indiana
Univ. Math. J., 38 (1989), pp. 861-915.

\bibitem{F-L-W}
F. Huang, M. Li, Y. Wang, Zero dissipation limit to rarefaction wave with
vacuum for one-dimensional compressible Navier-Stokes equations,
SIAM J. Math. Anal., 44 (2012),
pp. 1742-1759.

\bibitem{H-J-W} F.M. Huang, S. Jiang and Y. Wang, Zero dissipation limit of full compressible Navier-Stokes equations with a Riemann initial data, Preprint.

\bibitem{HWWY} F. Huang, Y. Wang, Y. Wang, T. Yang, The Limit of the Boltzmann Equation to the Euler Equations for Riemann Problems, SIAM J. Math. Anal., 45 (2013), no. 3, pp.1741-1811.



\bibitem{Huang-Wang-Yang}
F. M. Huang, Y. Wang and T. Yang, Fluid Dynamic Limit to the Riemann
Solutions of Euler Equations: I. Superposition of rarefaction waves
and contact discontinuity, Kinetic and Related Models, 3 (2010), pp.
685-728.

\bibitem{Huang-Wang-Yang-1}
F. M. Huang, Y. Wang and T. Yang, Vanishing viscosity limit of the
compressible Navier-Stokes equations for  solutions to Riemann
problem, Arch. Ration. Mech. Anal., 203 (2012),
pp. 379-413.




\bibitem{jiang}
S. Jiang, G. X. Ni and W. J. Sun, Vanishing viscosity limit to
rarefaction waves for the Navier-Stokes equations of one-dimensional
compressible heat-conducting fluids, SIAM J. Math. Anal., 38 (2006),
pp. 368-384.

\bibitem{JWX}Q. S. Jiu, Y. Wang and Z. P. Xin,
Vacuum behaviors around rarefaction waves to 1D compressible
Navier-Stokes equations with density-dependent viscosity, Preprint,
2011.

\bibitem{KMN} S. Kawashima, A. Matsumura, T. Nishida, On the fluid-dynamical approximation to the Boltzmann equation at the level of the Navier-Stokes equation, Comm. Math. Phys., 70 (1979), no. 2, pp. 97-124.

\bibitem{lw}
M. Li, T. Wang, Zero dissipation limit to rarefaction wave with vacuum for
 one-dimensional full compressible Navier-Stokes equations, Preprint, 2013.


\bibitem{lius}
T. P. Liu, J. Smoller,  On the vacuum state for the isentropic
 gas dynamics equations, Adv. in Appl. Math.,
 1 (1980), pp. 345-359.

 \bibitem{liuyyz}
 T. Liu, T. Yang, S. H. Yu and H. J. Zhao, Nonlinear stability of rarefaction waves for the Boltzmann
equation, Arch. Rat. Mech. Anal., 181 (2006), 333--371.

 \bibitem{Ma} S. X. Ma, Zero dissipation limit to strong
contact discontinuity for the 1-D compressible Navier-Stokes
equations, J. Differential Equations,  248 (2010), pp. 95--110.

\bibitem{mn86}
A. Matsumura, K. Nishihara, Asymptotics toward the rarefaction waves
of the solutions of a one-dimensional model system for compressible
viscous gas, Japan J. Appl. Math., 3 (1986), pp. 1-13.

\bibitem{P} M. Perepelitsa,  Asymptotics toward rarefaction waves and
vacuum for 1-D compressible Navier-Stokes equations, SIAM J. Math.
Anal., 42 (2010), pp. 1404-1412.

\bibitem{smoller}
J. Smoller,  Shock Waves and Reaction-Diffusion Equations. 2nd ed.
Grundlehren der Mathematischen Wissenschaften. 258. New York:
Springer-Verlag, xxii, 1994.

\bibitem{Wang-H} H. Y. Wang, Viscous limits for piecewise smooth solutions of the p-system, J. Math. Anal. Appl.,
299 (2004), pp. 411-432.

\bibitem{Wang} Y. Wang, Zero dissipation limit of the compressible heat-conducting Navier-Stokes equations in the presence of the shock,
Acta Mathematica Scientia, 28B (2008), pp. 727-748.

\bibitem{xin93}
Z. P. Xin, Zero dissipation limit to rarefaction waves for the
one-dimensional Navier-Stokes equations of compressible isentropic
gases, Comm. Pure Appl. Math., 46 (1993), pp. 621-665.

\bibitem{Xin-Zeng} Z. P. Xin, H. H. Zeng, Convergence to the rarefaction
waves for the nonlinear Boltzmann equation and compressible
Navier-Stokes equations, J. Diff. Eqs., 249 (2010), pp. 827-871.

\bibitem{Yu} S. H. Yu, Zero-dissipation limit of solutions with shocks
for systems of hyperbolic conservation laws, Arch. Ration. Mech.
Anal., 146 (1999), pp. 275-370.


\bibitem{ZPWT} Y. Zhang, R. Pan, Y. Wang, Z. Tan, Zero dissipation limit with two interacting shocks of the 1D
non-isentropic Navier-Stokes equations, to appear in Indiana Univ. Math. J..

\end{thebibliography}
\end{document}